\documentclass[12pt]{amsart}
\usepackage{amssymb}
\newtheorem{theorem}{Theorem}

\newtheorem{proposition}[theorem]{Proposition}
\newtheorem{conjecture}[theorem]{Conjecture}
\newtheorem{corollary}[theorem]{Corollary}
\newtheorem{remark}[theorem]{Remark}
\newtheorem{definition}{Definition}
\newtheorem{lemma}[theorem]{Lemma}

\newcommand{\Z}{{\mathbb Z}}

\newcommand{\R}{{\mathbb R}}
\newcommand{\C}{{\mathbb C}}

\linethickness{0.4pt}

\begin{document}


\title[Complements of discriminants of singularities $X_{10}$]{Complements of discriminants of real function singularities of type $X_{10}$}
\author{V.A.~Vassiliev}
\address{Weizmann Institute of Science}
\email{vavassiliev@gmail.com}
\date{}
\begin{abstract}
Conjecturally complete lists of
connected components of the complements of the discriminant varieties (aka wave fronts) of smooth function singularities of types $X_{10}^3$ and $X_{10}^1$ are provided. The combinatorial invariants of these components are defined and calculated. 

It is shown that the one-dimensional homology group of the complement of the discriminant
of any versal deformation of any singularity of class $X_{10}^1$ contains a subgroup isomorphic to $\Z^3$. 
\end{abstract}
\maketitle 

\section{Introduction}

\subsection{}
The {\em discriminant} of a family of smooth functions depending on parameters is the set of parameter values for which the corresponding functions have critical points with zero critical value, see e.g. \cite{AGLV}, \cite{AVGZ}.

Enumerating the connected components of the complements of 
discriminants of deformations of smooth function singularities is a local version of the classical problem of real algebraic geometry on rigid isotopy classification of smooth real hypersurfaces of a given degree, which comes back to Hilbert's 16th problem, see e.g. \cite{Kharlamov}, \cite{nik}, \cite{pol}, \cite{viro}.

Discriminants also appear as {\em wave fronts} in the theory of hyperbolic PDE's, see \cite{Petrovskii}, \cite{ABG}, \cite{Gording}, \cite{ramif}. The asymptotic behavior of the corresponding {\em wave functions}, i.e. the fundamental solutions of hyperbolic equations, is different in different components of the complement of the wave front. Therefore, the enumeration of these components is the first step in studying the local properties of wave functions.

In \cite{sedykh}, \cite{VLoo}, and \cite{parab}--\cite{parab2} this enumeration was performed for the deformations of the simplest function singularities, so-called {\em simple} and {\em parabolic} ones. Here, we present a conjectural solution to the same problem for the next in difficulty classes $X_{10}^1$ and $X_{10}^3$, providing lists of 44 and 14 components, respectively. We describe an algorithm for enumerating possible sets of topological invariants that characterize such components. For each set of invariants found by this algorithm, we construct perturbations of singularities of these classes that
have these invariants. The remaining problem is whether these invariants separate all possible components.

We also show that the complements of the discriminant varieties of the $X_{10}^1$ singularities have nontrivial one-dimensional homology groups, see Theorem \ref{onehom} on p.~\pageref{onehom}. 

Additionally, we find the first example in which the Gusein-Zade--A'Campo method (\cite{AC}, \cite{GZ}) of computing the Dynkin diagrams of plane function singularities yields a diagram with multiple edges, see Fig.~\ref{twenty} on page \pageref{twenty}.

\subsection{Main objects}
A {\em real function singularity} is a germ of a regular analytic function $(\C^n, \R^n, a) \to (\C, \R, 0)$ at some point $a \in \R^n$ such that $d f(a)=0$. Two function singularities at points $a$ and $b$ are {\em $\mbox{RL}$-equivalent} if they can be transformed into each other via a pair of local diffeomorphisms $(\R^n, a) \to (\R^n,b)$ 
and $(\R, 0) \to (\R,0)$ of their source and target spaces. 
Two function singularities, $f(x_1, \dots, x_n) $ and $\tilde f (\tilde x_1, \dots, \tilde x_m)$, are {\em stably $\mbox{\rm RL}$-equivalent} if they become equivalent after summing with non-degenerate quadratic forms in additional variables, that is, $$f(x_1, \dots, x_n) + Q_2(y_1, \dots, y_{N-n}) \sim \tilde f(\tilde x_1, \dots, \tilde x_m) + \tilde Q_2(\tilde y_1, \dots, \tilde y_{N-m})$$ 
for some $N$ and real quadratic polynomials $Q_2,$ $\tilde Q_2$ of maximal ranks.

Let $f:(\C^n,\R^n,0) \to (\C,\R,0)$ be a function germ with a critical point at the origin, and $F:(\C^n\times \C^k, \R^n \times \R^k,0) \to (\C, \R, 0)$ be some its {\em deformation}, i.e., a family of functions $f_\lambda\equiv F(\cdot, \lambda)$, $\lambda \in \C^k,$ defined in a neighborhood of the origin in $\C^n \times \C^k$, such that $f_0 \equiv f$ and $f_\lambda (\R^n) \subset \R$ for $\lambda \in \R^k$. The {\em real discriminant} $\Sigma(F) \subset \R^k$ of this deformation is the set of values $\lambda \in \R^k$ such that the function $f_\lambda$ has a critical point with zero critical value near the origin in $\R^n$. We study the local (near the origin in $\R^k$) connected components of the sets $\R^k \setminus \Sigma(F)$ for function singularities of classes $X_{10}^3$ and $X_{10}^1$ of the V.~Arnold's classification list. 

We will only consider sufficiently large, so-called {\em versal}, deformations of isolated function singularities (see, for example, \cite{AVGZ}), since all other deformations can be reduced to them. For all versal deformations of a given singularity, which depend on the same number $k$ of parameters, the germs of pairs $(\R^k, \Sigma(F))$ at the origin are locally diffeomorphic to each other. For versal deformations $F$ and $\tilde F$ depending respectively on $k$ and $k+l$ parameters the germ of the pair $(\R^{k+l}, \Sigma(\tilde F))$ is diffeomorphic to that of $(\R^k, \Sigma(F)) \times \R^l$.
Moreover, the ambient diffeomorphism type of the germs of the real discriminants of all $k$-parametric versal deformations of a function singularity depends only on the stable $RL$-equivalence class of the singularity. 

Therefore, it is sufficient to enumerate the local connected components of
the spaces $\R^k \setminus \Sigma(F)$ for only one function singularity from each stable $\mbox{RL}$-equivalence class, and for an arbitrary single versal deformation of this singularity.

\subsection{Previous results}

The natural classification of function singularities begins with the so-called {\em simple} and {\em parabolic} classes, see \cite{AVGZ}. The classes $X_{10}^3$ and $X_{10}^1$, which are studied below, are the next step in this classification. They are
the only two remaining singularity classes stably equivalent to functions in two variables with Milnor number at most 10. 

For all classes of codimension up to five in the functional space (i.e., classes $A_k$, $D_4$, $D_5$, $D_6$ and $E_6$, all of which are simple and stably equivalent to functions in two variables) the numbers of discriminant complements follow from the results of \cite{sedykh}. For all simple singularities these components were explicitly listed in \cite{VLoo}. This list repeats the list \cite{viro} of {\it topological types} of non-singular zero-level sets of perturbed functions of these classes. 

\begin{definition} \rm
\label{deftt}
Two small non-discriminant perturbations of the same isolated function singularity in two variables, defined in a neighborhood of the critical point, have the same {\em topological type} \label{toptype} if there is an isotopy between their zero-level sets in this neighborhood that keeps the behavior of these sets at its boundary. 
\end{definition}

The topological type is an invariant of the local connected components of the discriminant complements: the isotopies of the zero-level sets are provided by any paths connecting the functions within the complement of the discriminant. 
Thus, the main result of \cite{VLoo} states that it is a complete invariant of components of discriminant complements of simple singularities. Note also the work \cite{Lo}, in which a one-to-one correspondence was proved between the sets of these components and of some algebraic objects related to eponymous Weyl groups.
The works \cite{sedykh}, \cite{sed24}, \cite{sed25}, \cite{sed26} also contain a wealth of information on the geometry and stratification of discriminants of simple function singularities.

The complete list of components of discriminant complements of parabolic singularity classes (i.e., classes $P_8^1$, $P_8^2$, $X_9^{\pm}$, $X_9^1$, $X_9^2$, $J_{10}^1,$ and $J_{10}^3$) was given in \cite{parab2}. The topological type of the zero-level set of a perturbation of a parabolic singularity is not a complete invariant of components. The numbers of these topological types for these singularities in two variables (which were found in \cite{viro}) and of components of discriminant complements are compared in Table \ref{t0}. 
\begin{table}[h]
\caption{Numbers of rigid isotopy and topological classes of non-discriminant perturbations of parabolic singularities in two variables}
\label{t0}
\begin{center}
\begin{tabular}{|l|c|c|c|c|c|}
\hline
Singularity class & $X_9^{\pm}$ & $X_9^1$ & $X_9^2$ & $J_{10}^1$ & $J_{10}^3$ \\
\hline
No. of components & 7 & 14 & 52 & 13 & 33 \\
\hline
No. of topological types & 6 & 12 & 50 & 10 & 31 \\
\hline
\end{tabular}
\end{center}
\end{table}
In the present work, we give a conjectural solution to the same problem for the next in difficulty classes of singularities, $X_{10}^1$, and $X_{10}^3$.
Since these singularities are stably equivalent to the functions in two variables, we assume below that $n=2$. 
\medskip

All simple and parabolic singularities are {\em quasihomogeneous} in appropriate local coordinates, which significantly simplifies their study. The $X_{10}$ case is the first example of plane function singularities that do not satisfy this condition.

\section{Main results}
\subsection{Components of discriminant complements of singularities of type $X_{10}^3$}
\label{par3}

All function singularities of class $X_{10}^3$ in two variables can be reduced to the normal form 
\begin{equation}
\label{meqx3}
-x^5 + x^2y^2 -y^4
\end{equation}
by the action of the group $\mbox{RL}$ of local changes of coordinates in the source and target spaces. We will also use the equivalent normal form 
\begin{equation}
\label{meqx3a}
(y^2-x^3)(x^2-y^2) .
\end{equation}
The local (at the origin) zero-level set of either function (\ref{meqx3}) or (\ref{meqx3a}) consists of two smooth branches and a curve with a semicubical cusp, see Fig.~\ref{zeroset}.

\begin{remark} \rm
In all pictures in this paper, the $x$ coordinate is constant on the vertical lines and grows to the right; the $y$ coordinate is constant on the horizontal lines. 
\end{remark}

To not have to worry about the choice of the versal deformation of this singularity, we will assume that it is sufficiently large. Namely it will have the form \begin{equation}
\label{vers}
f_0 + \sum_{m+n \leq D} \lambda_{m,n} x^m y^n,
\end{equation} where $f_0$ is the initial function (\ref{meqx3}) or (\ref{meqx3a}), and $D$ is a large natural number.

\begin{figure}
\begin{picture}(50,30)
\put(2,0){\line(1,0){48}}
\put(2,0){\line(0,1){30}}
\put(50,30){\line(-1,0){48}}
\put(50,30){\line(0,-1){30}}
\put(2,0){\line(1,1){30}}
\put(2,30){\line(1,-1){30}}
\bezier{200}(50,27)(30,15)(17,15)
\bezier{200}(50,3)(30,15)(17,15)
\put(4,14){\tiny $+$}
\put(46,14){\tiny $-$}
\put(15,1){\tiny $-$}
\put(15,27){\tiny $-$}
\put(46,2){\makebox(0,0)[cc]{\tiny $+$}}
\put(46,27){\makebox(0,0)[cc]{\tiny $+$}}
\end{picture}
\caption{Zero set of a singularity of type $X_{10}^3$}
\label{zeroset}
\end{figure}

\begin{definition} \rm
\label{defind}
The {\it index} \ $\mbox{Ind}(f_\lambda)$ \ of a real non-discriminant strict morsification $f_\lambda$ of a function singularity is equal to the number of its critical points with negative critical value and even Morse index minus the number of critical points with negative critical value and odd Morse index. 
\end{definition}

\begin{proposition}
The numbers $\mbox{Ind}$ are invariants of components of the discriminant complement. 
\end{proposition}

\noindent
{\it Proof.} All standard surgeries of Morse functions that occur at the isotopies 
not crossing the discriminant preserve this number. \hfill $\Box$
\medskip

This invariant is determined by the topological type invariant, see Definition \pageref{deftt}. However, it is more convenient for the formal computations.

\begin{theorem}
\label{totalnum0}
The singularity $X_{10}^3$ has at least 44 different local components of the discriminant complement. The shapes of the local zero-level sets of the perturbations of the singularity (\ref{meqx3}) that
represent some 32 of these components are shown in Figs.~\ref{x3beg}--\ref{x3end}. 
The numbers in the right parts of the subscripts under  these pictures are equal to the $\mbox{Ind}$ invariants of corresponding components \footnote{The other, left-hand, parts of the subscripts present another invariant, which will be defined in \S \ref{basinv}.}. The shapes of the remaining twelve components can be obtained by the reflection across a horizontal line from the pictures with subscripts $(53232, \ 3)$, $(7200, \ 3)$, $(40608, \ 2)$, $(107662, \ 1)$, $(189496, \ 0)$, $(276604, \ 0)$, $(23756, \ -1)$, $(12848, \ -1)$, $(120088, \ -2)$, $(63700, \ -3)$, $(25104, \ -4)$, and $(4800, \ -5)$. 
\end{theorem}

The perturbations of singularity (\ref{meqx3}) with such zero-level sets will be constructed in \S \ref{real3}. 

\begin{proposition} 
All real morsifications representing the 44 pictures described in Theorem \ref{totalnum0} belong to different components of the complement of the discriminant.
\end{proposition}

{\it Proof.} For almost all pairs of these 44 pictures this follows from their topological non-equivalence. For the remaining three pairs: $(8496, \ -4)$ and $(9984, \ -4)$ in one case, $(81848, \ -4)$ and $(29040, \ -4)$ in the other case, and being the picture $(53232, \ 3)$ and its mirror image in the third case, this follows from Bezout's theorem. Indeed, 
the polynomial (\ref{meqx3}) has a versal deformation, all whose functions $f_\lambda$ are polynomials of degree 5. Therefore, for example, the real morsification whose shown in the picture $(53232, \ 3)$ cannot be continuously moved within this deformation to the morsification obtained from it by the involution $y \leftrightarrow -y$: any such movement of plane curves contains curves with more than five points on the same affine line. \hfill $\Box$

\begin{conjecture}
\label{mcj1}
The singularity $X_{10}^3$ has exactly 44 local connected components of the complement of the discriminant variety.
\end{conjecture}

\begin{figure}
\begin{center}
\unitlength 1mm
\begin{picture}(45.00,40.00)
\put(1,28){\small $+$}
\put(41,19){\small $-$}
\put(0,5){\line(1,0){45}}
\put(0,5){\line(0,1){30}}
\put(0,35){\line(1,0){45}}
\put(45,5){\line(0,1){30}}
\bezier{200}(1,35)(16,15)(31,35)
\bezier{200}(1,5)(16,25)(31,5)
\bezier{200}(45,30)(22,20)(45,10)
\put(8,20){\circle{4}}
\put(8,20){\small \makebox(0,0)[cc]{\small $-$}}
\put(17,20){\circle{4}}
\put(17,20){\small \makebox(0,0)[cc]{\small $-$}}
\put(26,16){\circle{4}}
\put(26,16){\small \makebox(0,0)[cc]{\small $-$}}
\put(15,32){\small $-$}
\put(15,6){\small $-$}
\put(1,1){$53232$}
\put(42,1){$3$}
\end{picture} \quad 
\begin{picture}(45.00,40.00)
\put(0.5,28){\small $+$}
\put(41,19){\small $-$}
\put(0,5){\line(1,0){45}}
\put(0,5){\line(0,1){30}}
\put(0,35){\line(1,0){45}}
\put(45,5){\line(0,1){30}}
\bezier{200}(45,28)(30,20)(23,18)
\bezier{200}(23,18)(15,16)(1,5)
\bezier{200}(1,35)(16,20)(31,35)
\bezier{200}(31,5)(25,20)(45,15)
\put(6.5,20){\circle{3}}
\put(6.5,20){\makebox(0,0)[cc]{\small $-$}}
\put(11,16){\circle{3}}
\put(11,16){\makebox(0,0)[cc]{\small $-$}}
\put(15,18){\circle{3}}
\put(15,18){\makebox(0,0)[cc]{\small $-$}}
\put(20,20){\circle{3}}
\put(20,20){\makebox(0,0)[cc]{\small $-$}}
\put(37,8){\small $+$}
\put(15,32){\small $-$}
\put(1,1){$7200$}
\put(42,1){$3$}
\end{picture} \quad
 \begin{picture}(45.00,40.00)
\put(1,28){\small $+$}
\put(41,19){\small $-$}
\put(0,5){\line(1,0){45}}
\put(0,5){\line(0,1){30}}
\put(0,35){\line(1,0){45}}
\put(45,5){\line(0,1){30}}
\bezier{200}(1,35)(16,15)(31,35)
\bezier{200}(1,5)(16,25)(31,5)
\bezier{200}(45,30)(22,20)(45,10)
\put(10,20){\circle{4}}
\put(10,20){\makebox(0,0)[cc]{\small $-$}}
\put(24,20){\circle{4}}
\put(24,20){\makebox(0,0)[cc]{\small $-$}}
\put(15,32){\small $-$}
\put(15,6){\small $-$}
\put(1,1){$214912$}
\put(42,1){$2$}
\end{picture} \quad 
\begin{picture}(45.00,40.00)
\put(1,28){\small $+$}
\put(41,19){\small $-$}
\put(0,5){\line(1,0){45}}
\put(0,5){\line(0,1){30}}
\put(0,35){\line(1,0){45}}
\put(45,5){\line(0,1){30}}
\bezier{200}(45,28)(30,20)(23,18)
\bezier{200}(23,18)(15,16)(1,5)
\bezier{200}(31,5)(25,20)(45,10)
\bezier{200}(1,35)(16,15)(31,35)
\put(8,20){\circle{3}}
\put(8,20){\makebox(0,0)[cc]{\small $-$}}
\put(17,20){\circle{3}}
\put(17,20){\makebox(0,0)[cc]{\small $-$}}
\put(24,23){\circle{3}}
\put(24,23){\makebox(0,0)[cc]{\small $-$}}
\put(37,8){\small $+$}
\put(15,32){\small $-$}
\put(1,1){$40608$}
\put(42,1){$2$}
\end{picture} \quad
\begin{picture}(45.00,40.00)
\put(1,28){\small $+$}
\put(41,19){\small $-$}
\put(0,5){\line(1,0){45}}
\put(0,5){\line(0,1){30}}
\put(0,35){\line(1,0){45}}
\put(45,5){\line(0,1){30}}
\bezier{200}(1,35)(16,15)(31,35)
\bezier{200}(1,5)(16,25)(31,5)
\bezier{200}(45,30)(22,20)(45,10)
\put(10,20){\circle{4}}
\put(10,20){\makebox(0,0)[cc]{\small $-$}}
\put(15,32){\small $-$}
\put(15,6){\small $-$}
\put(1,1){$385722$}
\put(42,1){$1$}
\end{picture} \quad 
\begin{picture}(45.00,40.00)
\put(1,28){\small $+$}
\put(41,19){\small $-$}
\put(0,5){\line(1,0){45}}
\put(0,5){\line(0,1){30}}
\put(0,35){\line(1,0){45}}
\put(45,5){\line(0,1){30}}
\bezier{200}(45,28)(30,20)(23,18)
\bezier{200}(23,18)(15,16)(1,5)
\bezier{200}(31,5)(25,20)(45,10)
\bezier{200}(1,35)(16,15)(31,35)
\put(9,20){\circle{4}}
\put(9,20){\makebox(0,0)[cc]{\small $-$}}
\put(24,23){\circle{4}}
\put(24,23){\makebox(0,0)[cc]{\small $-$}}
\put(37,8){\small $+$}
\put(15,32){\small $-$}
\put(1,1){$107662$}
\put(42,1){$1$}
\end{picture}
\begin{picture}(45.00,40.00)
\put(1,28){\small $+$}
\put(41,19){\small $-$}
\put(0,5){\line(1,0){45}}
\put(0,5){\line(0,1){30}}
\put(0,35){\line(1,0){45}}
\put(45,5){\line(0,1){30}}
\bezier{200}(1,35)(16,15)(31,35)
\bezier{200}(1,5)(16,25)(31,5)
\bezier{200}(45,30)(22,20)(45,10)
\put(15,32){\small $-$}
\put(15,6){\small $-$}
\put(1,1){$495364$}
\put(42,1){$0$}
\end{picture} \quad
\begin{picture}(45.00,40.00)
\put(1,28){\small $+$}
\put(41,19){\small $-$}
\put(0,5){\line(1,0){45}}
\put(0,5){\line(0,1){30}}
\put(0,35){\line(1,0){45}}
\put(45,5){\line(0,1){30}}
\bezier{200}(0,6)(15,20)(0,34)
\bezier{200}(31,35)(3,20)(31,5)
\bezier{200}(45,30)(22,20)(45,10)
\put(24,20){\circle{4}}
\put(24,20){\makebox(0,0)[cc]{\small $-$}}
\put(37,30){\small $+$}
\put(15,6){\small $-$}
\put(1,1){$7776$}
\put(42,1){$0$}
\end{picture} \quad
\begin{picture}(45.00,40.00)
\put(1,28){\small $+$}
\put(41,19){\small $-$}
\put(0,5){\line(1,0){45}}
\put(0,5){\line(0,1){30}}
\put(0,35){\line(1,0){45}}
\put(45,5){\line(0,1){30}}
\bezier{200}(1,5)(23,20)(1,35)
\bezier{200}(31,5)(15,20)(31,35)
\bezier{200}(45,10)(20,20)(45,30)
\put(7,20){\circle{4}}
\put(7,20){\makebox(0,0)[cc]{\small $-$}}
\put(37,30){\small $+$}
\put(15,6){\small $-$}
\put(1,1){$10784$}
\put(42,1){$0$}
\end{picture} \quad 
\begin{picture}(45.00,40.00)
\put(1,27){\small $+$}
\put(41,19){\small $-$}
\put(0,5){\line(1,0){45}}
\put(0,5){\line(0,1){30}}
\put(0,35){\line(1,0){45}}
\put(45,5){\line(0,1){30}}
\bezier{200}(45,28)(30,20)(23,18)
\bezier{200}(23,18)(15,16)(1,5)
\bezier{200}(31,5)(25,20)(45,10)
\bezier{200}(1,35)(16,15)(31,35)
\put(8,20){\circle{4}}
\put(8,20){\makebox(0,0)[cc]{\small $-$}}
\put(37,8){\small $+$}
\put(15,32){\small $-$}
\put(1,1){$189496$}
\put(42,1){$0$}
\end{picture} \quad
\begin{picture}(45.00,40.00)
\put(1,27){\small $+$}
\put(41,19){\small $-$}
\put(0,5){\line(1,0){45}}
\put(0,5){\line(0,1){30}}
\put(0,35){\line(1,0){45}}
\put(45,5){\line(0,1){30}}
\bezier{200}(0,34)(35,20)(0,6)
\bezier{200}(31,5)(15,23))(45,12)
\bezier{200}(31,35)(15,17))(45,28)
\put(6,20){\circle{4}}
\put(6,20){\makebox(0,0)[cc]{\small $-$}}
\put(12,20){\circle{4}}
\put(12,20){\makebox(0,0)[cc]{\small $-$}}
\put(37,8){\small $+$}
\put(37,30){\small $+$}
\put(1,1){$8976$}
\put(42,1){$0$}
\end{picture} 
\quad
\begin{picture}(45.00,40.00)
\put(1,28){\small $+$}
\put(41,19){\small $-$}
\put(0,5){\line(1,0){45}}
\put(0,5){\line(0,1){30}}
\put(0,35){\line(1,0){45}}
\put(45,5){\line(0,1){30}}
\bezier{200}(45,28)(30,20)(23,18)
\bezier{200}(23,18)(15,16)(1,5)
\bezier{200}(31,5)(25,20)(45,10)
\bezier{200}(1,35)(16,15)(31,35)
\put(15,32){\small $-$}
\put(37,8){\small $+$}
\put(1,1){$276604$}
\put(40,1){$-1$}
\end{picture}
\caption{Non-discriminant perturbations of the $X_{10}^3$ singularity (I)}
\label{x3beg}
\end{center}
\end{figure}
\begin{figure}
\begin{center}
\unitlength 1mm
 \begin{picture}(45.00,40.00)
\put(1,28){\small $+$}
\put(41,19){\small $-$}
\put(0,5){\line(1,0){45}}
\put(0,5){\line(0,1){30}}
\put(0,35){\line(1,0){45}}
\put(45,5){\line(0,1){30}}
\bezier{200}(1,35)(16,10)(31,35)
\bezier{200}(1,5)(16,20)(31,5)
\bezier{200}(45,10)(20,20)(45,30)
\put(16,28){\circle{5}}
\put(16,28){\makebox(0,0)[cc]{\small $+$}}
\put(22,33){\small $-$}
\put(15,6){\small $-$}
\put(1,1){$23756$}
\put(40,1){$-1$}
\end{picture} \quad 
\begin{picture}(45.00,40.00)
\put(1,28){\small $+$}
\put(41,19){\small $-$}
\put(0,5){\line(1,0){45}}
\put(0,5){\line(0,1){30}}
\put(0,35){\line(1,0){45}}
\put(45,5){\line(0,1){30}}
\bezier{200}(0,6)(15,20)(0,34)
\bezier{200}(31,35)(10,20)(31,5)
\bezier{200}(45,30)(22,20)(45,10)
\put(15,32){\small $-$}
\put(37,30){\small $+$}
\put(1,1){$123976$}
\put(40,1){$-1$}
\end{picture} \quad
\begin{picture}(45.00,40.00)
\put(1,28){\small $+$}
\put(41,19){\small $-$}
\put(0,5){\line(1,0){45}}
\put(0,5){\line(0,1){30}}
\put(0,35){\line(1,0){45}}
\put(45,5){\line(0,1){30}}
\bezier{200}(0,34)(35,20)(0,6)
\bezier{200}(31,5)(15,23))(45,12)
\bezier{200}(31,35)(15,17))(45,28)
\put(9,20){\circle{4}}
\put(9,20){\makebox(0,0)[cc]{\small $-$}}
\put(37,8){\small $+$}
\put(37,30){\small $+$}
\put(1,1){$92544$}
\put(40,1){$-1$}
\end{picture} 
 \begin{picture}(45.00,40.00)
\put(1,28){\small $+$}
\put(41,19){\small $-$}
\put(0,5){\line(1,0){45}}
\put(0,5){\line(0,1){30}}
\put(0,35){\line(1,0){45}}
\put(45,5){\line(0,1){30}}
\bezier{200}(45,28)(30,20)(23,18)
\bezier{200}(23,18)(15,16)(1,5)
\bezier{200}(31,5)(25,20)(45,10)
\bezier{200}(1,35)(16,15)(31,35)
\put(12,20){\circle{5}}
\put(12,20){\makebox(0,0)[cc]{\small $-$}}
\put(24,14){\circle{5}}
\put(24,14){\makebox(0,0)[cc]{\small $+$}}
\put(35,8){\small $+$}
\put(15,32){\small $-$}
\put(1,1){$12848$}
\put(40,1){$-1$}
\end{picture} \quad
\begin{picture}(45.00,40.00)
\put(1,28){\small $+$}
\put(41,19){\small $-$}
\put(0,5){\line(1,0){45}}
\put(0,5){\line(0,1){30}}
\put(0,35){\line(1,0){45}}
\put(45,5){\line(0,1){30}}
\bezier{200}(0,6)(15,20)(0,34)
\bezier{200}(31,35)(10,20)(31,5)
\bezier{200}(45,30)(22,20)(45,10)
\put(14,20){\circle{5}}
\put(14,20){\makebox(0,0)[cc]{\small $+$}}
\put(15,32){\small $-$}
\put(38,30){\small $+$}
\put(1,1){$108196$}
\put(40,1){$-2$}
\end{picture} \quad 
\begin{picture}(45.00,40.00)
\put(1,28){\small $+$}
\put(41,19){\small $-$}
\put(0,5){\line(1,0){45}}
\put(0,5){\line(0,1){30}}
\put(0,35){\line(1,0){45}}
\put(45,5){\line(0,1){30}}
\bezier{200}(45,28)(30,20)(23,18)
\bezier{200}(23,18)(15,16)(1,5)
\bezier{200}(31,5)(25,20)(45,10)
\bezier{200}(1,35)(16,15)(31,35)
\put(22,13){\circle{5}}
\put(22,13){\makebox(0,0)[cc]{\small $+$}}
\put(36,8){\small $+$}
\put(15,32){\small $-$}
\put(1,1){$120088$}
\put(40,1){$-2$}
\end{picture} 
 \begin{picture}(45.00,40.00)
\put(1,28){\small $+$}
\put(41,19){\small $-$}
\put(0,5){\line(1,0){45}}
\put(0,5){\line(0,1){30}}
\put(0,35){\line(1,0){45}}
\put(45,5){\line(0,1){30}}
\bezier{200}(0,34)(15,20)(0,6)
\bezier{200}(31,5)(15,23))(45,12)
\bezier{200}(31,35)(15,17))(45,28)
\put(37,8){\small $+$}
\put(37,30){\small $+$}
\put(1,1){$495992$}
\put(40,1){$-2$}
\end{picture} \quad 
\begin{picture}(45.00,40.00)
\put(1,28){\small $+$}
\put(41,19){\small $-$}
\put(0,5){\line(1,0){45}}
\put(0,5){\line(0,1){30}}
\put(0,35){\line(1,0){45}}
\put(45,5){\line(0,1){30}}
\bezier{200}(0,34)(30,20)(0,6)
\bezier{200}(31,5)(15,23))(45,12)
\bezier{200}(31,35)(15,17))(45,28)
\put(9,20){\circle{5}}
\put(9,20){\makebox(0,0)[cc]{\small $-$}}
\put(20,20){\circle{5}}
\put(20,20){\makebox(0,0)[cc]{\small $+$}}
\put(37,8){\small $+$}
\put(37,30){\small $+$}
\put(1,1){$55440$}
\put(40,1){$-2$}
\end{picture} \quad
\begin{picture}(45.00,40.00)
\put(1,28){\small $+$}
\put(41,19){\small $-$}
\put(0,5){\line(1,0){45}}
\put(0,5){\line(0,1){30}}
\put(0,35){\line(1,0){45}}
\put(45,5){\line(0,1){30}}
\bezier{200}(0,6)(15,20)(0,34)
\bezier{200}(31,35)(10,20)(31,5)
\bezier{200}(45,30)(22,20)(45,10)
\put(15,15){\circle{5}}
\put(15,15){\makebox(0,0)[cc]{\small $+$}}
\put(15,25){\circle{5}}
\put(15,25){\makebox(0,0)[cc]{\small $+$}}
\put(15,32){\small $-$}
\put(37,30){\small $+$}
\put(1,1){$67624$}
\put(40,1){$-3$}
\label{67624}
\end{picture} \qquad 
\begin{picture}(45.00,40.00)
\put(1,28){\small $+$}
\put(41,19){\small $-$}
\put(0,5){\line(1,0){45}}
\put(0,5){\line(0,1){30}}
\put(0,35){\line(1,0){45}}
\put(45,5){\line(0,1){30}}
\bezier{200}(45,28)(30,22)(23,20)
\bezier{200}(23,20)(15,18)(1,5)
\bezier{200}(1,35)(16,15)(31,35)
\bezier{200}(31,5)(15,20))(45,12)
\put(13,11){\circle{4}}
\put(13,11){\makebox(0,0)[cc]{\small $+$}}
\put(18,13){\circle{4}}
\put(18,13){\makebox(0,0)[cc]{\small $+$}}
\put(37,8){\small $+$}
\put(15,32){\small $-$}
\put(1,1){$63700$}
\put(40,1){$-3$}
\end{picture} \quad
\begin{picture}(45.00,40.00)
\put(1,28){\small $+$}
\put(41,19){\small $-$}
\put(0,5){\line(1,0){45}}
\put(0,5){\line(0,1){30}}
\put(0,35){\line(1,0){45}}
\put(45,5){\line(0,1){30}}
\bezier{200}(0,34)(20,20)(0,6)
\bezier{200}(31,5)(15,23))(45,12)
\bezier{200}(31,35)(15,17))(45,28)
\put(20,20){\circle{5}}
\put(20,20){\makebox(0,0)[cc]{\small $+$}}
\put(37,8){\small $+$}
\put(37,30){\small $+$}
\put(1,1){$303168$}
\put(40,1){$-3$}
\end{picture} \quad 
\begin{picture}(45.00,40.00)
\put(0.5,27){\small $+$}
\put(41,19){\small $-$}
\put(0,5){\line(1,0){45}}
\put(0,5){\line(0,1){30}}
\put(0,35){\line(1,0){45}}
\put(45,5){\line(0,1){30}}
\bezier{200}(0,34)(26,20)(0,6)
\bezier{200}(31,5)(15,17)(45,10)
\bezier{200}(31,35)(15,23)(45,30)
\put(20,17){\circle{4}}
\put(20,17){\makebox(0,0)[cc]{\small $+$}}
\put(20,23){\circle{4}}
\put(20,23){\makebox(0,0)[cc]{\small $+$}}
\put(5.8,20){\circle{4}}
\put(5.8,20){\makebox(0,0)[cc]{\small $-$}}
\put(35,8){\small $+$}
\put(35,30){\small $+$}
\put(1,1){$22800$}
\put(40,1){$-3$}
\end{picture} 
\begin{picture}(45.00,40.00)
\put(1,28){\small $+$}
\put(41,19){\small $-$}
\put(0,5){\line(1,0){45}}
\put(0,5){\line(0,1){30}}
\put(0,35){\line(1,0){45}}
\put(45,5){\line(0,1){30}}
\bezier{200}(0,6)(15,20)(0,34)
\bezier{200}(31,35)(10,20)(31,5)
\bezier{200}(45,30)(22,20)(45,10)
\put(17,20){\circle{4}}
\put(17,20){\makebox(0,0)[cc]{\small $+$}}
\put(13,13){\circle{4}}
\put(13,13){\makebox(0,0)[cc]{\small $+$}}
\put(13,27){\circle{4}}
\put(13,27){\makebox(0,0)[cc]{\small $+$}}
\put(15,32){\small $-$}
\put(37,30){\small $+$}
\put(1,1){$8496$}
\put(40,1){$-4$}
\label{8496p}
\end{picture} \quad 
\begin{picture}(45.00,40.00)
\put(1,28){\small $+$}
\put(41,19){\small $-$}
\put(0,5){\line(1,0){45}}
\put(0,5){\line(0,1){30}}
\put(0,35){\line(1,0){45}}
\put(45,5){\line(0,1){30}}
\bezier{200}(0,6)(15,20)(0,34)
\bezier{200}(31,35)(10,20)(31,5)
\bezier{200}(45,30)(22,20)(45,10)
\put(13,20){\circle{4}}
\put(13,20){\makebox(0,0)[cc]{\small $+$}}
\put(17,13){\circle{4}}
\put(17,13){\makebox(0,0)[cc]{\small $+$}}
\put(17,27){\circle{4}}
\put(17,27){\makebox(0,0)[cc]{\small $+$}}
\put(15,32){\small $-$}
\put(37,30){\small $+$}
\put(1,1){$9984$}
\put(40,1){$-4$}
\end{picture} \quad
\begin{picture}(45.00,40.00)
\put(1,28){\small $+$}
\put(41,19){\small $-$}
\put(0,5){\line(1,0){45}}
\put(0,5){\line(0,1){30}}
\put(0,35){\line(1,0){45}}
\put(45,5){\line(0,1){30}}
\bezier{200}(45,28)(30,22)(23,20)
\bezier{200}(23,20)(15,18)(1,5)
\bezier{200}(1,35)(16,15)(31,35)
\bezier{200}(31,5)(18,18))(45,12)
\put(13,11){\circle{4}}
\put(13,11){\makebox(0,0)[cc]{\small $+$}}
\put(18,14.5){\circle{4}}
\put(18,14.5){\makebox(0,0)[cc]{\small $+$}}
\put(23,15){\circle{4}}
\put(23,15){\makebox(0,0)[cc]{\small $+$}}
\put(35,8){\small $+$}
\put(15,32){\small $-$}
\put(1,1){$25104$}
\put(40,1){$-4$}
\end{picture} 
\caption{Non-discriminant perturbations of the $X_{10}^3$ singularity (II)}
\label{x3contin}
\end{center}
\end{figure}
\begin{figure}
\unitlength 1mm
\begin{center}
\begin{picture}(45.00,40.00)
\put(1,27){\small $+$}
\put(41,19){\small $-$}
\put(0,5){\line(1,0){45}}
\put(0,5){\line(0,1){30}}
\put(0,35){\line(1,0){45}}
\put(45,5){\line(0,1){30}}
\bezier{200}(0,34)(15,20)(0,6)
\bezier{200}(31,5)(15,23))(45,12)
\bezier{200}(31,35)(15,17))(45,28)
\put(17,13){\circle{5}}
\put(17,13){\makebox(0,0)[cc]{\small $+$}}
\put(17,27){\circle{5}}
\put(15.4, 26){\small $+$}
\put(37,8){\small $+$}
\put(37,30){\small $+$}
\put(1,1){$81848$}
\put(40,1){$-4$}
\end{picture} \quad
\begin{picture}(45.00,40.00)
\put(1,27){\small $+$}
\put(41,19){\small $-$}
\put(0,5){\line(1,0){45}}
\put(0,5){\line(0,1){30}}
\put(0,35){\line(1,0){45}}
\put(45,5){\line(0,1){30}}
\bezier{200}(0,34)(15,20)(0,6)
\bezier{200}(31,5)(15,23))(45,12)
\bezier{200}(31,35)(15,17))(45,28)
\put(13,20){\circle{5}}
\put(13,20){\makebox(0,0)[cc]{\small $+$}}
\put(20,20){\circle{5}}
\put(20,20){\makebox(0,0)[cc]{\small $+$}}
\put(37,8){\small $+$}
\put(37,30){\small $+$}
\put(1,1){$29040$}
\put(40,1){$-4$}
\end{picture} \quad
\begin{picture}(45.00,40.00)
\put(0.5,28){\small $+$}
\put(41,19){\small $-$}
\put(0,5){\line(1,0){45}}
\put(0,5){\line(0,1){30}}
\put(0,35){\line(1,0){45}}
\put(45,5){\line(0,1){30}}
\bezier{200}(0,34)(21,20)(0,6)
\bezier{200}(33,5)(18,23))(45,12)
\bezier{200}(33,35)(18,17))(45,28)
\put(19,13){\circle{4}}
\put(19,13){\makebox(0,0)[cc]{\small $+$}}
\put(19,27){\circle{4}}
\put(19, 27){\makebox(0,0)[cc]{\small $+$}}
\put(15,20){\circle{4}}
\put(15,20){\makebox(0,0)[cc]{\small $+$}}
\put(7,20){\circle{4}}
\put(7,20){\makebox(0,0)[cc]{\small $-$}}
\put(37,8){\small $+$}
\put(37,30){\small $+$}
\put(1,1){$4320$}
\put(40,1){$-4$}
\end{picture}
\begin{picture}(45.00,40.00)
\put(1,28){\small $+$}
\put(41,19){\small $-$}
\put(0,5){\line(1,0){45}}
\put(0,5){\line(0,1){30}}
\put(0,35){\line(1,0){45}}
\put(45,5){\line(0,1){30}}
\bezier{200}(45,28)(30,22)(23,20)
\bezier{200}(23,20)(15,18)(1,5)
\bezier{200}(1,35)(16,15)(31,35)
\bezier{200}(31,5)(18,20))(45,12)
\put(13,11){\circle{3}}
\put(13,11){\makebox(0,0)[cc]{\small $+$}}
\put(18,13){\circle{3}}
\put(18,13){\makebox(0,0)[cc]{\small $+$}}
\put(23,15){\circle{3}}
\put(23,15){\makebox(0,0)[cc]{\small $+$}}
\put(8,8){\circle{3}}
\put(8,8){\makebox(0,0)[cc]{\small $+$}}
\put(37,8){\small $+$}
\put(15,32){\small $-$}
\put(1,1){$4800$}
\put(40,1){$-5$}
\end{picture} \quad 
\begin{picture}(45.00,40.00)
\put(1,28){\small $+$}
\put(41,19){\small $-$}
\put(0,5){\line(1,0){45}}
\put(0,5){\line(0,1){30}}
\put(0,35){\line(1,0){45}}
\put(45,5){\line(0,1){30}}
\bezier{200}(0,34)(15,20)(0,6)
\bezier{200}(31,5)(15,23))(45,12)
\bezier{200}(31,35)(15,17))(45,28)
\put(13,20){\circle{4}}
\put(13,20){\makebox(0,0)[cc]{\small $+$}}
\put(18,13){\circle{4}}
\put(18,13){\makebox(0,0)[cc]{\small $+$}}
\put(18,27){\circle{4}}
\put(18,27){\makebox(0,0)[cc]{\small $+$}}
\put(37,8){\small $+$}
\put(37,30){\small $+$}
\put(15,32){\small $-$}
\put(15,6){\small $-$}
\put(1,1){$13392$}
\put(40,1){$-5$}
\end{picture} 
\caption{Non-discriminant perturbations of the $X_{10}^3$ singularity (III)}
\label{x3end}
\end{center}
\end{figure}

\subsection{Answers for singularities of class $X_{10}^1$}

All singularities of class $X_{10}^1$ are $\mbox{RL}$-equivalent to the function 
\begin{equation}
(x^2+y^2)(x^3-y^2) .
\label{meqx101}
\end{equation}
The zero-level set of this function in $\R^2$ is the semicubical parabola $y^2 = x^3$. 

\begin{theorem}
\label{mthm1010}
There are at least 14 distinct local connected components of the discriminant complement of the $X_{10}^1$ singularity. The shapes of the zero-level sets of perturbations of $($\ref{meqx101}$)$ representing some twelve of them are given in Fig. \ref{tta}. The shapes of such sets of remaining two perturbations can be obtained by the reflection in a horizontal line from the pictures 
$(12672, \ 4)$ and $(84936, \ 3)$.
\end{theorem}

\begin{figure}
{\rm \begin{center}
\unitlength 1mm
\begin{picture}(45.00,40.00)
\put(1,6){\small $-$}
\put(41,29){\small $+$}
\put(0,5){\line(1,0){45}}
\put(0,5){\line(0,1){30}}
\put(0,35){\line(1,0){45}}
\put(45,5){\line(0,1){30}}
\bezier{350}(25,13)(-20,20)(25,27)
\bezier{150}(45,5)(38,11)(25,13)
\bezier{150}(45,35)(38,29)(25,27)
\put(12,19){\makebox(0,0)[cc]{\small $-$}}
\put(22,22){\makebox(0,0)[cc]{\small $-$}}
\put(32,16){\makebox(0,0)[cc]{\small $-$}}
\put(42,20){\makebox(0,0)[cc]{\small $-$}}
\put(12,19){\circle{4}}
\put(22,22){\circle{4}}
\put(32,16){\circle{4}}
\put(42,20){\circle{4}}
\put(1,1){12672}
\put(42,1){ 4}
\end{picture} \quad 
\begin{picture}(45.00,40.00)
\put(1,6){\small $-$}
\put(41,30){\small $+$}
\put(0,5){\line(1,0){45}}
\put(0,5){\line(0,1){30}}
\put(0,35){\line(1,0){45}}
\put(45,5){\line(0,1){30}}
\bezier{400}(45,34)(-25,20)(45,6)
\put(19,19){\makebox(0,0)[cc]{\small $-$}}
\put(28,22){\makebox(0,0)[cc]{\small $-$}}
\put(40,16){\makebox(0,0)[cc]{\small $-$}}
\put(19,19){\circle{5}}
\put(28,22){\circle{5}}
\put(40,16){\circle{5}}
\put(1,1){84936}
\put(42,1){3}
\end{picture} \quad
\begin{picture}(45.00,40.00)
\put(1,6){\small $-$}
\put(41,26){\small $+$}
\put(0,5){\line(1,0){45}}
\put(0,5){\line(0,1){30}}
\put(0,35){\line(1,0){45}}
\put(45,5){\line(0,1){30}}
\bezier{400}(45,31)(-25,20)(45,9)
\put(20,20){\makebox(0,0)[cc]{\small $-$}}
\put(32,20){\makebox(0,0)[cc]{\small $-$}}
\put(20,20){\circle{5}}
\put(32,20){\circle{5}}
\put(1,1){202812}
\put(42,1){2}
\end{picture} 
\begin{picture}(45.00,40.00)
\put(1,6){\small $-$}
\put(41,26){\small $+$}
\put(0,5){\line(1,0){45}}
\put(0,5){\line(0,1){30}}
\put(0,35){\line(1,0){45}}
\put(45,5){\line(0,1){30}}
\bezier{400}(45,31)(-25,20)(45,9)
\put(20,20){\makebox(0,0)[cc]{\small $-$}}
\put(20,20){\circle{5}}
\put(1,1){277887}
\put(42,1){1} 
\end{picture} \quad
\begin{picture}(45.00,40.00)
\put(1,6){\small $-$}
\put(41,26){\small $+$}
\put(0,5){\line(1,0){45}}
\put(0,5){\line(0,1){30}}
\put(0,35){\line(1,0){45}}
\put(45,5){\line(0,1){30}}
\bezier{400}(45,31)(-25,20)(45,9)
\put(1,1){526228}
\put(42,1){0}
\end{picture} \quad 
\begin{picture}(45.00,40.00)
\put(1,6){\small $-$}
\put(41,26){\small $+$}
\put(0,5){\line(1,0){45}}
\put(0,5){\line(0,1){30}}
\put(0,35){\line(1,0){45}}
\put(45,5){\line(0,1){30}}
\bezier{400}(45,31)(-5,20)(45,9)
\put(35,20){\circle{5}}
\put(35,20){\makebox(0,0)[cc]{\small $-$}}
\put(10,20){\circle{5}}
\put(10,20){\makebox(0,0)[cc]{\small $+$}}
\put(1,1){7776}
\put(42,1){0}
\end{picture} 
\begin{picture}(45.00,40.00)
\put(1,6){\small $-$}
\put(41,26){\small $+$}
\put(0,5){\line(1,0){45}}
\put(0,5){\line(0,1){30}}
\put(0,35){\line(1,0){45}}
\put(45,5){\line(0,1){30}}
\bezier{400}(45,31)(15,20)(45,9)
\put(16,20){\makebox(0,0)[cc]{\small $-$}}
\put(16,24.5){\makebox(0,0)[cc]{\small $+$}}
\put(16,20){\circle{5}}
\put(16,20){\circle{15}}
\put(1,1){7536}
\put(42,1){0}
\end{picture}
\quad 
\begin{picture}(45.00,40.00)
\put(1,6){\small $-$}
\put(41,26){\small $+$}
\put(0,5){\line(1,0){45}}
\put(0,5){\line(0,1){30}}
\put(0,35){\line(1,0){45}}
\put(45,5){\line(0,1){30}}
\bezier{400}(45,31)(-5,20)(45,9)
\put(10,20){\makebox(0,0)[cc]{\small $+$}}
\put(10,20){\circle{5}}
\put(1,1){306917}
\put(40,1){$-1$}
\end{picture} \quad 
\begin{picture}(45.00,40.00)
\put(1,6){\small $-$}
\put(41,26){\small $+$}
\put(0,5){\line(1,0){45}}
\put(0,5){\line(0,1){30}}
\put(0,35){\line(1,0){45}}
\put(45,5){\line(0,1){30}}
\bezier{400}(45,31)(15,20)(45,9)
\put(10,20){\makebox(0,0)[cc]{\small $+$}}
\put(20,20){\makebox(0,0)[cc]{\small $+$}}
\put(10,20){\circle{5}}
\put(20,20){\circle{5}}
\put(1,1){265228}
\put(40,1){$-2$}
\end{picture} 
\begin{picture}(45.00,40.00)
\put(1,6){\small $-$}
\put(41,26){\small $+$}
\put(0,5){\line(1,0){45}}
\put(0,5){\line(0,1){30}}
\put(0,35){\line(1,0){45}}
\put(45,5){\line(0,1){30}}
\bezier{400}(45,31)(15,20)(45,9)
\put(10,20){\makebox(0,0)[cc]{\small $+$}}
\put(20,25){\makebox(0,0)[cc]{\small $+$}}
\put(20,15){\makebox(0,0)[cc]{\small $+$}}
\put(10,20){\circle{5}}
\put(20,25){\circle{5}}
\put(20,15){\circle{5}}
\put(1,1){157004}
\put(40,1){$-3$ }
\end{picture} \quad
\begin{picture}(45.00,40.00)
\put(1,6){\small $-$}
\put(41,26){\small $+$}
\put(0,5){\line(1,0){45}}
\put(0,5){\line(0,1){30}}
\put(0,35){\line(1,0){45}}
\put(45,5){\line(0,1){30}}
\bezier{400}(45,31)(15,20)(45,9)
\put(16,20){\makebox(0,0)[cc]{\small $+$}}
\put(24,28){\makebox(0,0)[cc]{\small $+$}}
\put(24,12){\makebox(0,0)[cc]{\small $+$}}
\put(7,20){\makebox(0,0)[cc]{\small $+$}}
\put(7,20){\circle{5}}
\put(16,20){\circle{5}}
\put(24,28){\circle{5}}
\put(24,12){\circle{5}}
\put(1,1){19296}
\put(40,1){$-4$}
\end{picture} \quad 
\begin{picture}(45.00,40.00)
\put(1,6){\small $-$}
\put(41,26){\small $+$}
\put(0,5){\line(1,0){45}}
\put(0,5){\line(0,1){30}}
\put(0,35){\line(1,0){45}}
\put(45,5){\line(0,1){30}}
\bezier{400}(45,31)(15,20)(45,9)
\put(16,28){\makebox(0,0)[cc]{\small $+$}}
\put(24,20){\makebox(0,0)[cc]{\small $+$}}
\put(16,12){\makebox(0,0)[cc]{\small $+$}}
\put(7,20){\makebox(0,0)[cc]{\small $+$}}
\put(7,20){\circle{5}}
\put(16,12){\circle{5}}
\put(24,20){\circle{5}}
\put(16,28){\circle{5}}
\put(1,1){20736}
\put(40,1){$-4$}
\end{picture}
\caption{Non-discriminant perturbations of the $X_{10}^1$ singularity}
\label{tta}
\end{center}
}
\end{figure}

All these twelve pictures represent functions from different components of the discriminant complement. Indeed, almost all of them are topologically distinct. The only two pictures with the same topological type (having $\mbox{Ind}=-4$) cannot be connected in the space of non-discriminant perturbations of (\ref{meqx101}) due to Bezout's theorem reasons. These reasons also prove that the first picture of Fig. \ref{tta} (with \ $\mbox{Ind}=4$) and its image under the involution $y \leftrightarrow -y$ belong to different components of the complement of the discriminant. The same is true also for the second picture (with \ $\mbox{Ind}=3$).

\begin{conjecture}
\label{mcj2}
The singularity $X_{10}^1$ has exactly 14 local connected components of the complement of the discriminant variety, which are mentioned in Theorem \ref{mthm1010}.
\end{conjecture}

\begin{theorem}
\label{onehom}
The one-dimensional homology group of the discriminant complement 
of $X_{10}^1$ singularity has a subgroup isomorphic to $\Z^3$. Specifically, each component containing the functions shown in pictures $(265228, -2)$, $(157004, -3),$ and $(19296, -4)$ has a summand of $\Z$ in its one-dimensional homology group.
\end{theorem}

\section{Realization of components for the $X_{10}^3$ singularity}
\label{real3}

\subsection{Realization by perturbed sabirizations}

\label{realsab}

\subsubsection{}
Let $f: (\C^2, \R^2, 0) \to (\C, \R,0)$ be an analytic function with a critical point at the origin, and $\mu(f)$ be its Milnor number. Let $\tilde f$ be a small strict morsification of $f$, all $\mu(f)$ critical points of which are real (and are distinct since the function is {\em strictly} Morse).

\begin{definition} 
\label{defscale}
\rm The {\em standard scale} of $\tilde f$ is any collection of $\mu(f)+1$ strictly Morse non-discriminant functions, all of which are equal to $\tilde f$ and to each other up to the addition of constant functions, chosen so that the numbers of negative critical values of these functions take all values from $0$ to $\mu(f)$.
\end{definition}

We will consider the scales of the strict morsifications that are the small perturbations of special (non-strict) morsifications of the following kind.

\begin{definition} \rm
\label{sabdef}
A {\em sabirization} of a real function singularity $f$ is any its small Morse perturbation which has exactly $\mu(f)$ real critical points, such that all its saddlepoints have critical value 0, and all its minima (respectively, maxima) have negative (respectively, positive) critical values. 
\end{definition}

Such morsifications are a crucial tool in the Gusein-Zade--A'Campo method for computing Dynkin diagrams of real function singularities, see \cite{AC}, \cite{GZ}.

\subsubsection{}
\label{ssseven}

\begin{figure}
\begin{picture}(50,40)
\put(0,5){\line(1,0){50}}
\put(0,5){\line(0,1){30}}
\put(50,35){\line(-1,0){50}}
\put(50,35){\line(0,-1){30}}
\bezier{300}(50,30)(20,0)(9,9)
\bezier{200}(9,9)(-5,20)(9,31)
\bezier{300}(50,10)(20,40)(9,31)
\put(2,5){\line(1,1){30}}
\put(2,35){\line(1,-1){30}}
\put(8,19){\tiny B}
\put(27,19){\tiny A}
\put(4.5,9.5){\tiny F}
\put(4,28.4){\tiny E}
\put(16,17.3){\tiny G}
\put(25.3,7.6){\tiny D}
\put(24.7,30.7){\tiny C}
\put(40.5,19){\tiny H}
\put(15.3,10.3){\tiny J}
\put(15.3,26.5){\tiny I}
\put(7.5,10.5){\circle*{1}}
\put(7.5,29.5){\circle*{1}}
\put(17,20){\circle*{1}}
\put(26.3,10.5){\circle*{1}}
\put(26.3,29.5){\circle*{1}}
\put(39.55,19.85){\circle*{1}}
\put(1,7){\tiny $+$}
\put(6,6){\tiny $-$}
\put(6,33){\tiny $-$}
\put(46,6){\tiny $+$}
\put(47,20){\tiny $-$}
\put(46, 33){\tiny $+$}
\end{picture} \qquad
\begin{picture}(60,36)
\put(0,0){\line(1,0){60}}
\put(0,0){\line(0,1){36}}
\put(60,36){\line(-1,0){60}}
\put(60,36){\line(0,-1){36}}
\bezier{400}(60,36)(7,8)(3,17)
\bezier{150}(3,17)(2.5,18)(3,19)
\bezier{400}(60,0)(7,28)(3,19)
\put(40,0){\line(1,2){18}}
\put(44,0){\line(-1,2){18}}
\put(9,17){\tiny A}
\put(29,17){\tiny B}
\put(41,6.5){\tiny C}
\put(20.7,15.3){\tiny G}
\put(32.2,24){\tiny E}
\put(39.9,3.3){\tiny I}
\put(39.2,11){\tiny F}
\put(57,31.5){\tiny D}
\put(45,8){\tiny H}
\put(43,20){\tiny J}
\put(21.5,18){\circle*{1}}
\put(32.9,22.6){\circle*{1}}
\put(38.8,10.8){\circle*{1}}
\put(44,8.2){\circle*{1}}
\put(57.5,34.5){\circle*{1}}
\put(42,4){\circle*{1}}
\put(1,1){\tiny $+$}
\put(52,1){\tiny $+$}
\put(57,17){\tiny $-$}
\put(46, 33){\tiny $-$}
\end{picture}
\caption{Two sabirizations of the $X_{10}^3$ singularity}
\label{general}
\end{figure}

It is easy to perturb the polynomial (\ref{meqx3a}) in such a way that the zero-level set of the resulting function near the origin looks as in Fig.~\ref{general} (left), where the letters $A$ and $B$ indicate the minima, the letters $I$ and $J$ indicate the maxima, and the remaining letters denote saddlepoints at self-crossings of the zero set. We can slightly modify this function so that the critical values of the new function $\tilde f$ at the new critical points (neighboring to the old ones) become ordered {\em alphabetically}: 
\begin{equation}
\label{standa}
\tilde f(\tilde A)<\tilde f(\tilde B)<\tilde f(\tilde C)<\tilde f(\tilde D)<\tilde f(\tilde E)<\tilde f(\tilde F)<\tilde f(\tilde G)<\tilde f(\tilde H)<\tilde f(\tilde I)<\tilde f(\tilde J). \end{equation}The standard scale of the obtained morsification realizes the pictures given in Figs.~\ref{x3beg}--\ref{x3end} with the subscripts {\bf (495364, \ 0), (385722, \ 1), (214912, \ 2)}, (385722, \ 1), {\bf (10784, \ 0), (123976, \ --1), (108196, \ --2), (67624, \ --3), (81848, \ --4), (303168, \ --3)}, and {\bf (495922, \ --2)}. (Here and below, we highlight in bold the titles of pictures that are being  realized for the first time.)

\subsubsection{}
\label{ssfour}

We can also perturb the function (\ref{meqx3a}) so that the topological picture of the zero-level set will be as shown in Fig.~\ref{general} (right), and then move it slightly so that the critical values are ordered according to the rule (\ref{standa}). The standard scale of the obtained perturbation realizes the pictures $(495364, \ 0)$, $(385722, \ 1)$, $(214912, \ 2)$, ${\bf (53232, \ 3)}$, ${\bf (40608, \ 2)},$ ${\bf (107662, \ 1)}$, ${\bf (189496, \ 0)}$, ${\bf (276604, \ -1)}$, ${\bf (120088, \ -2)}$, $(303168, \ -3)$, and $(495992, \ -2)$.

\subsubsection{} 
\label{sseight}

Let us take the same perturbation of the zero-level set as in Fig.~\ref{general} (right), but move the critical values so that 
$\tilde f(\tilde A)<\tilde f(\tilde B)<\tilde f(\tilde C)<\tilde f(\tilde D)<\tilde f(\tilde H)<\tilde f(\tilde I)<\tilde f(\tilde G)<\tilde f(\tilde F)<\tilde f(\tilde E)<\tilde f(\tilde J).$
The corresponding standard scale contains the functions shown in the pictures (495364, \ 0), (385722, \ 1), (214912, \ 2), (53232, \ 3), (40608, \ 2), (107662, \ 1), {\bf (8976, \ 0), (92544, \ --1)}, (495992, \ --2), (303168, \ --3), and (495992, \ --2). 

\subsubsection{}
\label{sssix}

Consider the perturbation of $f$, whose zero-level set is shown in Fig.~\ref{six} (left), and move its critical values so that they are ordered according to (\ref{standa}). 
The corresponding standard scale then consists of the functions shown in the pictures (495364, \ 0), (385722, \ 1), (214912, \ 2), (53232, \ 3), (214912, \ 2), (385722, \ 1), (495364, \ 0), {\bf (23756, \ --1),} (108196, \ --2), (303168, \ --3), and (495992, \ --2).

\begin{figure}
\unitlength 1.1mm
\begin{picture}(60,45)
\put(0,0){\line(1,0){60}}
\put(0,0){\line(0,1){42}}
\put(60,42){\line(-1,0){60}}
\put(60,42){\line(0,-1){42}}
\bezier{400}(60,41)(30,6)(10,18)
\bezier{250}(10,18)(0,26)(10,34)
\bezier{400}(10,34)(30,46)(60,11)
\put(20,0){\line(1,2){21}}
\put(26,0){\line(-1,2){21}}
\put(8,23){\tiny A}
\put(21.7,10.6){\tiny B}
\put(36,24){\tiny C}
\put(9,35){\tiny D}
\put(16.5,12.2){\tiny E}
\put(28,13){\tiny F}
\put(33.7,33.8){\tiny G}
\put(44.5,27){\tiny I}
\put(19.7,5.1){\tiny H}
\put(22,26){\tiny J}
\put(9.3,33.1){\circle*{1}}
\put(18.6,15){\circle*{1}}
\put(23,6){\circle*{1}}
\put(28,16){\circle*{1}}
\put(36.2,32.3){\circle*{1}}
\put(45,26){\circle*{1}}
\put(1,1){\tiny $+$}
\put(56,1){\tiny $+$}
\put(55,26){\tiny $-$}
\put(21,39){\tiny $-$}
\put(49,38){\tiny $+$}
\end{picture}
\qquad
\begin{picture}(60,45)
\put(0,0){\line(1,0){60}}
\put(0,0){\line(0,1){39}}
\put(60,39){\line(-1,0){60}}
\put(60,39){\line(0,-1){39}}
\bezier{300}(60,38)(20,14)(10,19)
\bezier{300}(60,10)(20,34)(10,29)
\bezier{150}(10,19)(0,24)(10,29)
\put(36,0){\line(-2,3){26}}
\put(30,0){\line(2,3){26}}
\put(32,14){\tiny A}
\put(11,22){\tiny B}
\put(51.8,30.2){\tiny C}
\put(33, 25){\tiny D}
\put(42.4,16.2){\tiny E}
\put(34,3){\tiny F}
\put(16,31){\tiny G}
\put(20.9,16.2){\tiny H}
\put(23,23){\tiny I}
\put(41,23){\tiny J}
\put(33,4.5){\circle*{1}}
\put(52,33){\circle*{1}}
\put(34,24){\circle*{1}}
\put(43,19.5){\circle*{1}}
\put(16,30){\circle*{1}}
\put(23,19.5){\circle*{1}}
\put(1,1){\tiny $+$}
\put(55,1){\tiny $+$}
\put(32,36){\tiny $-$}
\put(55,23){\tiny $-$}
\end{picture}
\caption{Two other sabirizations of the $X_{10}^3$ singularity}
\label{six}
\end{figure}

\subsubsection{}

\label{sstwo}

Take the perturbation of $f$, whose zero-level set is shown in Fig.~\ref{six} (right), and move the critical values so that they become ordered alphabetically.
The standard scale of the resulting morsification provides the pictures (495364, \ 0), (385722, \ 1), (214912, \ 2), (107662, \ 1), (189496, \ 0), {\bf (12848, \ --1), (55440, \ --2)}, (303168, \ --3), {\bf (29040, \ --4)}, (303168, \ --3), and (495992, \ --2).

\subsubsection{}
\label{sstop}

\begin{figure}
\begin{picture}(60,50)
\put(0,0){\line(1,0){60}}
\put(0,0){\line(0,1){40}}
\put(60,40){\line(-1,0){60}}
\put(60,40){\line(0,-1){40}}
\bezier{400}(60,5)(35,50)(25,30)
\bezier{400}(60,35)(35,-10)(25,10)
\bezier{200}(25,30)(18,20)(25,10)
\put(0,10){\line(3,2){45}}
\put(0,30){\line(3,-2){45}}
\put(33,20){\tiny A}
\put(14,17){\tiny B}
\put(19.8,24.9){\tiny C}
\put(19.8,12.8){\tiny D}
\put(36,36){\tiny E}
\put(36,1.9){\tiny F}
\put(52,19){\tiny G}
\put(18.5,18.7){\tiny H}
\put(29,6){\tiny J}
\put(30,32){\tiny I}
\put(15,20){\circle*{1}}
\put(22.5,15){\circle*{1}}
\put(22.5,25){\circle*{1}}
\put(37.2,5.2){\circle*{1}}
\put(37.2,34.8){\circle*{1}}
\put(51,20){\circle*{1}}
\put(1,1){\tiny $-$}
\put(55,1){\tiny $+$}
\put(1,37){\tiny $-$}
\put(57,20){\tiny $-$}
\put(1,20){\tiny $+$}
\put(55,37){\tiny $+$}
\end{picture}
\caption{The fifth sabirization of $X_{10}^3$}
\label{top}
\end{figure}

Consider the perturbation of $f$, whose zero-level set is shown in Fig.~\ref{top} and again arrange the critical values in alphabetical order. The standard scale of the obtained perturbation provides the pictures (495364, \ 0), (385722, \ 1), {\bf (7776, \ 0)}, (123976, \ --1), (108196, \ --2), (67624, \ --3), {\bf (9984, \ --4), (13392, \ --5)}, (81848, \ --4), (303168, \ --3), and (495992, \ --2).

\subsection{Realization by adjacent singularities}
\label{realadj}

\subsubsection{$X_{10}^3 \to E_8$}
\label{adje8}

The small perturbation $f \mapsto f+ \varepsilon y^3 ,$ where $\varepsilon >0$, moves the polynomial (\ref{meqx3}) to a function with a critical point of class $E_8$ at the origin. The zero-level set of this function near the origin looks as shown in Fig.~\ref{toe83} (right top), with the $E_8$ critical point marked by a thick dot. 
\begin{figure}
\begin{center}
\unitlength 1.2 mm
\begin{picture}(58,30)
\put(2,0){\line(1,0){48}}
\put(2,0){\line(0,1){30}}
\put(50,30){\line(-1,0){48}}
\put(50,30){\line(0,-1){30}}
\put(2,0){\line(1,1){30}}
\put(2,30){\line(1,-1){30}}
\bezier{200}(50,27)(30,15)(17,15)
\bezier{200}(50,3)(30,15)(17,15)
\put(4,14){\small $+$}
\put(46,14){\small $-$}
\put(15,1){\small $-$}
\put(15,27){\small $-$}
\put(44,3){\makebox(0,0)[cc]{\small $+$}}
\put(44,27){\makebox(0,0)[cc]{\small $+$}}
\put(52,15){$\longrightarrow$}
\end{picture} 
\unitlength 1.12 mm
\begin{picture}(45.00,36.00)
\put(0,0){\line(1,0){45}}
\put(0,0){\line(0,1){32}}
\put(45,32){\line(-1,0){45}}
\put(45,32){\line(0,-1){32}}
\bezier{300}(15,12)(30,12)(45,24)
\bezier{300}(15,12)(5,12)(0,0)
\bezier{300}(0,32)(15,8)(30,32)
\bezier{300}(30,0)(20,16)(45,6)
\put(15,12){\circle*{1}}
\put(15,1){\small $-$}
\put(15,29){\small $-$}
\put(36,2){\small $+$}
\put(1,24){\small $+$}
\put(40,28){\small $+$}
\end{picture}

\bigskip

\begin{picture}(45.00,43.00)
\put(1,28){\small $+$}
\put(41,19){\small $-$}
\put(14,6){\small $-$}
\put(14,32){\small $-$}
\put(35,6){\small $+$}
\put(35,32){\small $+$}
\put(7,19){\small $-$}
\put(0,5){\line(1,0){45}}
\put(0,5){\line(0,1){30}}
\put(0,35){\line(1,0){45}}
\put(45,5){\line(0,1){30}} 
\put(1,5){\line(1,1){30}}
\put(1,35){\line(1,-1){30}}
\bezier{300}(16,20)(25,23)(45,33)
\bezier{300}(16,20)(25,17)(45,7)
\bezier{150}(16,20)(10,17)(6,18)
\bezier{120}(6,18)(2,20)(6,22)
\bezier{150}(16,20)(10,23)(6,22)
\put(22,40){$\downarrow$}
\put(53,40){$\searrow$}
\end{picture} \qquad \qquad
\begin{picture}(45.00,35.00)
\put(0,5){\line(1,0){45}}
\put(0,5){\line(0,1){30}}
\put(0,35){\line(1,0){45}}
\put(45,5){\line(0,1){30}} 
\put(1,5){\line(1,1){30}}
\put(1,35){\line(1,-1){30}}
\bezier{250}(45,10)(15, 20)(45,30)
\put(2,19){\small $+$}
\put(40,19){\small $-$}
\put(13,6){\small $-$}
\put(13, 32){\small $-$}
\put(22,19){\small $+$}
\end{picture}
\caption{Perturbations $X_{10}^3 \to E_8$ (top), 
$X_{10}^3 \to X_9^2$ (left bottom) and $X_{10}^3 \to X_9^1$ (right bottom)} 
\label{toe83}
\end{center}
\end{figure}
This its critical point has ten topologically distinct further perturbations, which are listed in Proposition 3.6 of \cite{VLoo}. Applying all of these perturbations yields, we obtain ten pictures of zero-level sets that are topologically equivalent to the pictures {\bf (7200, \ 3)}, (40608, \ 2), (107662, \ 1), (189496, \ 0), (276604, \ --1), (12848, \ --1), (120088, \ --2), {\bf (63700, \ --3)}, {\bf (25104, \ --4)}, and {\bf (4800, \ --5)} on Figs.~\ref{x3beg}--\ref{x3end}.

\subsubsection{$X_{10}^3 \to X_9^2$}
\label{tox92}

The perturbation $\tilde f \equiv (y^2-x^3- \varepsilon x^2)(x^2-y^2),$ $\varepsilon >0$, of the polynomial (\ref{meqx3a}) has a critical point of class $X_9^2$ at the origin. The topology of the zero-level set of this perturbation $\tilde f$ is shown in Fig.~\ref{toe83}. In particular, the perturbation $\tilde f$ has one Morse minimum point in addition to the singularity at the origin. 

According to \cite{viro} and \cite{parab}, any  singularity of class $X_9^2$ (whose zero-level set is shown in the center of Fig.~\ref{pertx92}) 
\begin{figure}
\unitlength 1.15mm
\begin{picture}(45.00,34.00)
\put(0,5){\line(1,0){30}}
\put(0,5){\line(0,1){30}}
\put(0,35){\line(1,0){30}}
\put(30,5){\line(0,1){30}}
\bezier{200}(30,25)(15,20)(20,35)
\bezier{200}(30,15)(15,20)(20,5)
\bezier{200}(10,5)(15,20)(10,35)
\bezier{200}(0,15)(15,20)(0,25)
\put(15,23){\circle{2.3}}
\put(15,17){\circle{2.3}}
\put(15,23){\makebox(0,0)[cc]{\tiny $+$}}
\put(15,17){\makebox(0,0)[cc]{\tiny $+$}}
\put(26,6){\small $+$}
\put(26,32){\small $+$}
\put(1,19){\small $-$}
\put(1,32){\small $+$}
\put(14,32){\small $-$}
\put(34,18){$\longleftarrow$}
\end{picture} 
\begin{picture}(45.00,34.00)
\put(0,5){\line(1,0){30}}
\put(0,5){\line(0,1){30}}
\put(0,35){\line(1,0){30}}
\put(30,5){\line(0,1){30}}
\put(10,5){\line(1,3){10}}
\put(20,5){\line(-1,3){10}}
\put(0,15){\line(3,1){30}}
\put(0,25){\line(3,-1){30}}
\put(1,32){\small $+$}
\put(1,6){\small $+$}
\put(26,32){\small $+$}
\put(26,6){\small $+$}
\put(14,6){\small $-$}
\put(14,32){\small $-$}
\put(1,19){\small $-$}
\put(26,19){\small $-$}
\put(34,18){$\longrightarrow$}
\end{picture} 
\begin{picture}(30.00,35.00)
\put(0,5){\line(1,0){30}}
\put(0,5){\line(0,1){30}}
\put(0,35){\line(1,0){30}}
\put(30,5){\line(0,1){30}}
\bezier{200}(30,26)(16,20)(21,35)
\bezier{200}(30,14)(16,20)(21,5)
\bezier{200}(8,5)(13,20)(8,35)
\bezier{200}(0,15)(13,20)(0,25)
\put(16,24){\circle{2.3}}
\put(14,20){\circle{2.3}}
\put(16,16){\circle{2.3}}
\put(16,24){\makebox(0,0)[cc]{\tiny $+$}}
\put(14,20){\makebox(0,0)[cc]{\tiny $+$}}
\put(16,16){\makebox(0,0)[cc]{\tiny $+$}}
\put(26,6){\small $+$}
\put(26,32){\small $+$}
\put(1,19){\small $-$}
\put(1,32){\small $+$}
\put(14,32){\small $-$}
\end{picture} 
\caption{Perturbations of $X_9^2$ singularities}
\label{pertx92}
\end{figure}
has a perturbation shown in the same figure on the left and also three other perturbations whose pictures are obtained from this one by rotations to the angles $\pi/2$, $\pi$, and $3\pi/2$. For instance, this picture can be obtained from the picture (3840, 1) in Fig.~5 of \cite{parab} by multiplying the function by $-1$. Applying these four perturbations to the critical point of $\tilde f$ shown in Fig.~\ref{toe83} (bottom left), we obtain the pictures (67624, \ --3), (63700, \ --3), {\bf (22800, \ --3)}, and the mirror image of the picture (63700, \ --3). 

Any $X_9^2$ singularity also has a perturbation shown in Fig.~\ref{pertx92} on the right and also three perturbations obtained by rotating it to the angles $\pi/2$, $\pi$ and $3\pi/2$. Namely, the zero-level set of a preliminary perturbation is shown in different scalings on Fig.~\ref{realX92C}, see Fig.~32 of \cite{vasX9}. Further, we can slightly move its critical values in such a way that all the critical values at the saddlepoints except for the leftmost one (marked by 9 in Fig.~\ref{realX92C}) will be negative, and the critical values at the local maxima inside the three bounded domains will remain positive. Moreover, we can arrange the order of all critical values as shown by the numbers in Fig.~\ref{realX92C}.
\unitlength 0.4mm
\begin{figure}
\begin{picture}(100,100)
\put(0,0){\line(1,0){100}}
\put(0,0){\line(0,1){100}}
\put(100,100){\line(-1,0){100}}
\put(100,100){\line(0,-1){100}}
\put(70,50){\line(1,3){16.7}}
\put(70,50){\line(1,-3){16.7}}
\put(70,50){\line(-1,3){16.7}}
\put(70,50){\line(-1,-3){16.7}}
\bezier{700}(100,68)(70.1,62)(70.1,50)
\bezier{700}(100,32)(70.1,38)(70.1,50)
\bezier{700}(0,65)(25,60)(25,50)
\bezier{700}(0,35)(25,40)(25,50)
\put(45,50){\circle*{2.5}}
\put(45,44){\makebox(0,0)[cc]{\tiny $9$}}
\put(5,50){\makebox(0,0)[cc]{\tiny $-$}}
\put(95,50){\makebox(0,0)[cc]{\tiny $-$}}
\put(70,4){\makebox(0,0)[cc]{\tiny $-$}}
\put(70,95){\makebox(0,0)[cc]{\tiny $-$}}
\put(5,5){\makebox(0,0)[cc]{\tiny $+$}}
\put(95,5){\makebox(0,0)[cc]{\tiny $+$}}
\put(95,95){\makebox(0,0)[cc]{\tiny $+$}}
\end{picture} \qquad
\begin{picture}(100,100)
\put(0,0){\line(1,0){100}}
\put(0,0){\line(0,1){100}}
\put(100,100){\line(-1,0){100}}
\put(100,100){\line(0,-1){100}}
\put(60,50){\line(1,3){16.7}}
\put(60,50){\line(1,-3){16.7}}
\put(60,50){\line(-1,3){16.7}}
\put(60,50){\line(-1,-3){16.7}}
\bezier{700}(100,100)(61,85)(61,50)
\bezier{700}(100,0)(61,15)(61,50)
\put(5,50){\circle*{2.5}}
\put(5,44){\makebox(0,0)[cc]{\tiny $9$}}
\put(69.5,78){\circle*{2.5}}
\put(73.5,77){\makebox(0,0)[cc]{\tiny $1$}}
\put(69.5,22){\circle*{2.5}}
\put(73.5,23){\makebox(0,0)[cc]{\tiny $5$}}
\put(88,60){\vector(-3,1){23}}
\put(91.5,60){\makebox(0,0)[cc]{\tiny $6$}}
\put(88,40){\vector(-3,-1){23}}
\put(91.5,40){\makebox(0,0)[cc]{\tiny $8$}}
\end{picture} \qquad
\begin{picture}(100,100)
\put(0,0){\line(1,0){100}}
\put(0,0){\line(0,1){100}}
\put(100,100){\line(-1,0){100}}
\put(100,100){\line(0,-1){100}}
\put(60,50){\line(1,3){16.7}}
\put(60,50){\line(1,-3){16.7}}
\put(60,50){\line(-1,3){16.7}}
\put(60,50){\line(-1,-3){16.7}}
\bezier{700}(70,100)(64,65)(64,50)
\bezier{700}(70,0)(64,35)(64,50)
\put(64.5,64){\circle*{2.5}}
\put(68.5,64){\makebox(0,0)[cc]{\tiny $2$}}
\put(64.5,36){\circle*{2.5}}
\put(60,50){\circle*{2.5}}
\put(68.5,36){\makebox(0,0)[cc]{\tiny $4$}}
\put(56,50){\makebox(0,0)[cc]{\tiny $3$}}
\put(83.3,50){\vector(-1,0){21}}
\put(86.3,50){\makebox(0,0)[cc]{\tiny $7$}}
\end{picture}
\caption{Realization of the right perturbation of Fig.~\ref{pertx92}}
\label{realX92C}
\end{figure}
The intersection matrix of the vanishing cycles of the standard basis corresponding to these critical points is as follows: 
\begin{equation} 
\begin{array}{ccccccccc}
 $-2$ & 0 & 0 & 0 & 0 & 1 & 0 & 0 & 0 \\
 0 & $-2$ & 1 & 0 & 0 & 1 & 1 & 0 & 0 \\
 0 & 1 & $-2$ & 1 & 0 & 0 & $-2$ & 0 & 1 \\
 0 & 0 & 1 & $-2$ & 0 & 0 & 1 & 1 & 0 \\
 0 & 0 & 0 & 0 & $-2$ & 0 & 0 & 1 & 0 \\
 1 & 1 & 0 & 0 & 0 & $-2$ & 0 & 0 & 0 \\
 0 & 1 & $-2$ & 1 & 0 & 0 & $-2$ & 0 & $1$ \\
 0 & 0 & 0 & 1 & 1 & 0 & 0 & $-2$ & 0 \\
 0 & 0 & 1 & 0 & 0 & 0 & 1 & 0 & $-2$ \\
\end{array} \ , \label{4320a}
\end{equation}
see the corresponding Dynkin graph in Fig.~16 of \cite{vasX9}.

By applying these four perturbations to the germ of $\tilde f$ at the origin, we obtain the following pictures: {$\bf (4320, \ -4)$}, $(25104, \ -4)$, {$\bf (8496, \ -4)$}, and the mirror image of the picture $(25104, \ -4)$ with respect to a horizontal line. 

All of 32 pictures from Figs.~\ref{x3beg}--\ref{x3end} are thus realized. Since the polynomial (\ref{meqx3}) is even in the $y$ coordinate, we can apply the perturbations symmetric to these 32 ones with respect to the horizontal axis to it. Twelve of the resulting components will be new; see the final sentence 
of Theorem \ref{totalnum0}. This theorem is thus proved. \hfill $\Box$

\subsubsection{$X_{10}^3 \to X_9^1$}
\label{tox91}

The perturbation $(y^2-x^3+ \varepsilon x^2)(x^2-y^2),$ $\varepsilon >0$, of the function (\ref{meqx3a}) has a critical point point of type $X_9^1$ at the origin, see Fig.~\ref{toe83} (right bottom).
This critical point has two further arbitrarily small perturbations, which are shown in the upper row of Fig.~\ref{x91pert}. The construction of the perturbation shown on the right is illustrated in the bottom part of the figure.
\begin{figure}
\begin{center}
\unitlength 1.8mm
\begin{picture}(30.00,27.00)
\put(0,5){\line(1,0){20}}
\put(0,5){\line(0,1){20}}
\put(0,25){\line(1,0){20}}
\put(20,5){\line(0,1){20}}
\bezier{200}(0,6)(10,15)(0,24)
\bezier{200}(20,6)(10,15)(20,24)
\put(12.5,8){\circle{2.5}}
\put(12.5,22){\circle{2.5}}
\put(9,15){\circle{2.5}}
\put(0.5,9){\scriptsize $+$}
\put(17.5,9){\scriptsize $+$}
\put(3,22){\scriptsize $-$}
\put(12.5,8){\makebox(0,0)[cc]{\tiny $+$}}
\put(12.5, 22){\makebox(0,0)[cc]{\tiny $+$}}
\put(9, 15){\makebox(0,0)[cc]{\tiny $+$}}
\put(24,13){$\longleftarrow$}
\end{picture} \
\begin{picture}(30.00,25.00)
\put(0,5){\line(1,0){20}}
\put(0,5){\line(0,1){20}}
\put(0,25){\line(1,0){20}}
\put(20,5){\line(0,1){20}}
\put(0,5){\line(1,1){20}}
\put(20,5){\line(-1,1){20}}
\put(0.5,8){\scriptsize $+$}
\put(18,8){\scriptsize $+$}
\put(10,24){\makebox(0,0)[cc]{\tiny $-$}}
\put(10,6){\makebox(0,0)[cc]{\tiny $-$}}
\put(24,13){$\longrightarrow$}
\end{picture}
\begin{picture}(20.00,25.00)
\put(0,5){\line(1,0){20}}
\put(0,5){\line(0,1){20}}
\put(0,25){\line(1,0){20}}
\put(20,5){\line(0,1){20}}
\bezier{200}(0,6)(10,15)(0,24)
\bezier{200}(20,6)(10,15)(20,24)
\put(8,8){\circle{2.5}}
\put(8,22){\circle{2.5}}
\put(11.5,15){\circle{2.5}}
\put(0.5,9){\scriptsize $+$}
\put(18,9){\scriptsize $+$}
\put(14,22){\scriptsize $-$}
\put(8,8){\makebox(0,0)[cc]{\tiny $+$}}
\put(8, 22){\makebox(0,0)[cc]{\tiny $+$}}
\put(11.5, 15){\makebox(0,0)[cc]{\tiny $+$}}
\end{picture}

\unitlength 2.8mm
\begin{picture}(20.00,28.00)
\put(9,14){\circle*{0.4}}
\put(9,16){\circle*{0.4}}
\put(6,11){\circle*{0.4}}
\put(6,19){\circle*{0.4}}
\put(10,15){\circle*{0.4}}
\put(5,15){\makebox(0,0)[cc]{\tiny $1$}}
\put(9.3,13.2){\makebox(0,0)[cc]{\tiny $6$}}
\put(11.3,15){\makebox(0,0)[cc]{\tiny $2$}}
\put(9.3,16.8){\makebox(0,0)[cc]{\tiny $3$}}
\put(7,20){\makebox(0,0)[cc]{\tiny $5$}}
\put(7,10){\makebox(0,0)[cc]{\tiny $4$}}
\put(11.7,21.7){\makebox(0,0)[cc]{\tiny $8$}}
\put(11.7,8.3){\makebox(0,0)[cc]{\tiny $7$}}
\put(13.7,19.7){\makebox(0,0)[cc]{\tiny $9$}}
\put(11,21){\vector(-1,-1){3.25}}
\put(11,9){\vector(-1,1){3.25}}
\put(13,19){\vector(-1,-1){3.75}}
\put(0,5){\line(1,0){20}}
\put(0,5){\line(0,1){20}}
\put(0,25){\line(1,0){20}}
\put(20,5){\line(0,1){20}}
\put(0,5){\line(1,1){20}}
\put(20,5){\line(-1,1){20}}
\put(0.5,8){\scriptsize $+$}
\put(17.5,8){\scriptsize $+$}
\put(3,24){\makebox(0,0)[cc]{\tiny $-$}}
\put(3,6){\makebox(0,0)[cc]{\tiny $-$}}
\bezier{200}(5,11)(1,11)(1,15)
\bezier{200}(1,15)(1,19)(5,19)
\bezier{200}(5,19)(9,19)(9,15)
\bezier{200}(9,15)(9,11)(5,11)

\put(10,28){{\large $\downarrow$}}
\put(23,27){$\nearrow$}
\end{picture}
\caption{Perturbations of the $X_9^1$ singularities}
\label{x91pert}
\end{center}
\end{figure}

After applying these perturbations, we  obtain pictures (9984, \ --4) and (8496, \ --4) of Fig.~\ref{x3contin} once again. 

\section{Realization of components for the $X_{10}^1$ singularity} 
\label{real4}

\begin{figure}
\unitlength 1.1mm
\begin{picture}(60,46)
\put(0,0){\line(1,0){60}}
\put(0,0){\line(0,1){36}}
\put(60,36){\line(-1,0){60}}
\put(60,36){\line(0,-1){36}}
\bezier{500}(60,36)(7,8)(3,17)
\bezier{150}(3,17)(2.5,18)(3,19)
\bezier{500}(60,0)(7,28)(3,19)
\bezier{200}(33,18)(33,33)(45,33)
\bezier{200}(45,33)(57,33)(57,18)
\bezier{200}(57,18)(57,3)(45,3)
\bezier{200}(45,3)(33,3)(33,18)
\put(7,16){\tiny J}
\put(28.3,16){\tiny G}
\put(19.6,13.7){\tiny B}
\put(41,28){\tiny I}
\put(39.9,4.4){\tiny H}
\put(30,9.7){\tiny E}
\put(53.2,29.6){\tiny C}
\put(54,3.1){\tiny D}
\put(30,23){\tiny F}
\put(42,17){\tiny A}
\put(50,17){\small $-$}
\put(21.5,18){\circle*{1}}
\put(33.4,22.8){\circle*{1}}
\put(33.4,13.2){\circle*{1}}
\put(51.3,4.3){\circle*{1}}
\put(51.3,31.5){\circle*{1}}
\put(5,5){\small $-$}
\end{picture} \qquad 
\begin{picture}(60,46)
\put(0,0){\line(1,0){60}}
\put(0,0){\line(0,1){36}}
\put(60,36){\line(-1,0){60}}
\put(60,36){\line(0,-1){36}}
\bezier{400}(60,33)(30,-2)(10,10)
\bezier{250}(10,10)(0,18)(10,26)
\bezier{400}(10,26)(30,38)(60,3)
\bezier{200}(12,18)(12,34)(24,34)
\bezier{200}(24,34)(36,34)(36,18)
\bezier{200}(36,18)(36,2)(24,2)
\bezier{200}(24,2)(12,2)(12,18)
\put(7,15.3){\tiny G}
\put(22.2,2.7){\tiny J}
\put(38,16.2){\tiny I}
\put(10.5,28.5){\tiny C}
\put(11,3.7){\tiny F}
\put(35.2,7.2){\tiny E}
\put(35.2,25.6){\tiny D}
\put(43.8,19.2){\tiny B}
\put(22.2,29.8){\tiny H}
\put(21.8,15.3){\tiny A}
\put(13.8,27.7){\circle*{1}}
\put(13.85,7.9){\circle*{1}}
\put(35.4,11.2){\circle*{1}}
\put(35.4,24.8){\circle*{1}}
\put(45,18){\circle*{1}}
\put(4,4){\small $-$}
\put(45,32){\small $-$}
\put(55,18){\small $+$}
\put(45,4){\small $-$}
\end{picture}
\caption{Two sabirizations of the $X_{10}^1$ singularity}
\label{fourteen}
\end{figure}

\begin{figure}
\unitlength 1.2mm
\begin{picture}(50,58)
\put(0,0){\line(1,0){50}}
\put(0,0){\line(0,1){50}}
\put(50,50){\line(-1,0){50}}
\put(50,50){\line(0,-1){50}}
\bezier{200}(5,27)(5,9)(23,9)
\bezier{200}(23,9)(41,9)(41,27)
\bezier{200}(41,27)(41,45)(23,45)
\bezier{200}(23,45)(5,45)(5,27)
\put(34.5,11.8){\circle*{1}}
\put(30,36.5){\circle*{1}}
\put(35.8,41.2){\circle*{1}}
\put(40.7,31.8){\circle*{1}}
\put(40.7,22.3){\circle*{1}}
\put(46,45){\small $+$}
\put(2,2){\small $-$}
\bezier{300}(45,50)(22,27)(15,32)
\bezier{150}(15,32)(13,34)(14,36)
\bezier{250}(14,36)(18,43)(38,33)
\bezier{250}(38,33)(60,24)(37,22)
\bezier{250}(37,22)(10,19)(34,12)
\bezier{250}(34,12)(45,8)(50,0)
\put(20,23){\tiny J}
\put(19, 34){\tiny C}
\put(35,36){\tiny D}
\put(43.5,25){\tiny B}
\put(34,9){\tiny F}
\put(41,20.5){\tiny G}
\put(41,32){\tiny E}
\put(34.5,42.5){\tiny H}
\put(29,34){\tiny I}
\put(31,15.5){\tiny A}
\put(10,16){\small $+$}
\end{picture} \quad 
\begin{picture}(50,58)
\put(0,0){\line(1,0){50}}
\put(0,0){\line(0,1){50}}
\put(50,50){\line(-1,0){50}}
\put(50,50){\line(0,-1){50}}
\bezier{200}(5,22)(5,4)(23,4)
\bezier{200}(23,4)(41,4)(41,22)
\bezier{200}(41,22)(41,40)(23,40)
\bezier{200}(23,40)(5,40)(5,22)
\bezier{200}(45,50)(39,39)(29,38.5)
\bezier{200}(32,38.4)(40,34)(50,16)
\bezier{60}(32,38.4)(30.2,39)(29,39)
\bezier{40}(29,38.5)(28,38.7)(29,39)
\put(35,40){\line(-1,0){10}}
\put(35,40){\line(0,-1){6}}
\put(25,34){\line(1,0){10}}
\put(25,34){\line(0,1){6}}
\end{picture} \quad
\begin{picture}(50,50)
\put(0,0){\line(1,0){50}}
\put(0,0){\line(0,1){40}}
\put(50,40){\line(-1,0){50}}
\put(50,40){\line(0,-1){40}}
\bezier{300}(50,40)(38,25.6)(31,27.2)
\bezier{300}(50,16)(38,29.6)(31,28.8)
\bezier{150}(31,27.2)(29,28)(31,28.8)
\bezier{150}(40,25.4)(35,29.4)(27,33.4)
\bezier{150}(27,33.4)(17,38.2)(0, 39.8)
\bezier{150}(40,25.6)(45,21.6)(50,13.6)
\put(33,20.8){\line(1,0){11}}
\put(33,20.8){\line(0,1){12.8}}
\put(44,33.6){\line(-1,0){11}}
\put(44,33.6){\line(0,-1){12.8}}
\end{picture} \quad
\begin{picture}(50,55)
\put(0,0){\line(1,0){50}}
\put(0,0){\line(0,1){50}}
\put(50,50){\line(-1,0){50}}
\put(50,50){\line(0,-1){50}}
\bezier{400}(50,50)(37,23)(0,20)
\bezier{400}(50,0)(37,27)(0,30)
\bezier{400}(0,37)(29,27)(50,2)
\end{picture}
\caption{Another sabirization of the $X_{10}^1$ singularity}
\label{twenty}
\end{figure}

We use three sabirizations of the singularity (\ref{meqx101}). All of them are products of two polynomials of degrees two and three, which are perturbations of two factors of the polynomial (\ref{meqx101}). Two of these sabirizations are shown in Fig.~\ref{fourteen}, and the third one is topologically situated as shown in Fig.~\ref{twenty} (left top) and consists of a degree-three curve  with a tiny loop, and a much larger circle. See Fig.~\ref{twenty} for its pictures at different scales, where the small rectangle in the top right picture corresponds to the entire bottom left picture, and the rectangle in the bottom left picture corresponds to the entire bottom right picture. 

We slightly perturb each of these three sabirizations in such a way that their critical values become ordered in the standard (alphabetical) way (\ref{standa}).

The shapes of the elements of the standard scale of the perturbed sabirization shown in Fig.~\ref{fourteen} $($left$)$ are given in the pictures $\bf (526228, 0)$, $\bf (277887, \ 1)$, $\bf (7776, \ 0)$, $\bf (306917, \ -1)$, $\bf (265228, -2)$, $\bf (157004, \ -3)$, $\bf (19296, \ -4)$, $(157004, \ -3)$, $(265228, \ -2)$, $(306917, \ -1)$, and $(526228, 0)$ of Fig.~\ref{tta}.
The analogous list for the standard scale of Fig.~\ref{fourteen} $($right$)$ consists of the pictures $ (526228, 0)$, $(277887, 1),$ $\bf (7536, \ 0)$, $(306917, \ -1)$, $(265228, \ -2)$, $(157004, \ -3)$, $\bf (20736, \ -4 )$, $(157004, \ -3)$, $(265228, \ -2)$, $(306917, \ -1)$, and $(526228, 0).$ 
The list for Fig.~\ref{twenty} consists of the pictures 
$(526228, 0)$, $(277887, \ 1),$ $\bf (202812, \ 2),$ $\bf (84936, \ 3)$, $\bf (12672, \ 4),$ $(84936, \ 3)$, $(202812, \ 2)$, $(277887, \ 1)$, $(7536, \ 0)$, $(306917, \ -1)$, and $(526228, 0)$.

In particular, all of the pictures in Fig.~\ref{tta} appear at these three scales. 
 Thus, Theorem \ref{mthm1010} is proved. \hfill $\Box$

\section{Proof of Theorem \ref{onehom}}

\label{onehompr}

 Recall the notation $B(\R^2, k)$ for the {\em unordered $k$-th configuration space} of the plane, whose points are arbitrary sets of $k$ pairwise distinct points in $\R^2$. For any $k \geq 2$, the group $H^1(B(\R^2,k))$ is isomorphic to $\Z$, see \cite{A70}.
The basic element of this group is induced from the generator of the group $H^1(\C^1 \setminus \{0\}) \sim \Z $ via the map $B(\R^2, k) \equiv B(\C^1,k) \to \C^1 \setminus \{0\}$ that sends any collection of \ $k$ \ distinct points $z_1, \dots, z_k \in \C^1$ to the product of all \ $k(k-1)$ \ complex numbers $z_i - z_j$, $i \neq j$. Let $U$ be a connected component of the discriminant complement of $X_{10}^1$ singularity that consists of functions whose zero-level sets include exactly $k \geq 2$ non-nested ovals, such that the functions are positive inside these ovals. (Such sets for different $k$ are shown in pictures $(265228, -2)$, $(157004, -3)$, $(19296, -4)$, and $(20736, -4)$ of Fig. \ref{tta}). Then, a map $H^*(B(\R^2,k)) \to H^*(U) $ is defined as follows. For any polynomial $f_\lambda \in U,$ define $\beta (f_\lambda)$ as the space of all choices of a point inside any of its $k$ ovals. These spaces $\beta(f_\lambda)$ form a locally trivial fiber bundle over the component $U$. The fibers of this fiber bundle are contractible; therefore, it has a cross-section, which is a continuous map $U \to B(\R^2, k)$. These maps defined by different cross-sections are clearly homotopic to each other and define the same homomorphism of the cohomology rings. 
Define the class $W_k \in H^1(U)$ as the cohomology class obtained by this map from the generator of the group $H^1(B(\R^2,k))$. 

\begin{proposition}
\label{prohom}
The cohomology class $W_2$ $($respectively, $W_3$, $W_4)$ takes all integer values on the appropriate 1-cycles in the component $(265228, -2)$ $($respectively, $(157004, -3)$, $(19296, -4))$ of the discriminant complements of the $X_{10}^1$ singularities.
\end{proposition}

\noindent
{\it Proof.} Consider the perturbation $f_0 \mapsto f_0 - \tau x^4$
 of the function (\ref{meqx101}) where $\tau$ is a small positive number. The perturbed function has a critical point of type \ $X_9^-$ \ at the origin. The space of functions (\ref{vers}) provides a versal deformation of the critical point of type $X_9^-$ of this perturbation, see \cite{parab2}.

In particular, the discriminant complement of a versal deformation of the $X_9^-$ singularity is embedded to the discriminant complement of the $X_{10}^1$ singularity. The analog of Proposition \ref{prohom} is proven in \cite{parab2} for some three components of the discriminant complement of any versal deformation of the $X_9^-$ singularity, see Theorem 3 of \cite{parab2}. The images of these three components under the embedding belong to the three components considered in Proposition \ref{prohom}. 
Thus, the integer cohomology classes $W_2, W_3, $ and $W_4$ take non-zero values on one-dimensional cycles contained in an arbitrary neighborhood of the function $f_0 - \tau x^4$ in the space of versal deformations of the class $X_{10}^1$ singularity $f_0$. \hfill $\Box$ 

\section{The search for components and their basic invariants}

In this section we describe our method for finding components of discriminant complements, cf. \cite{VLoo}, \cite{parab2}. We also define invariants of these components and calculate their values for all the components in the lists of Theorems \ref{totalnum0} and \ref{mthm1010}.

\subsection{The idea}
Our method consists in 
studying a graph, whose vertices are the collections of topological characteristics of
possible topological types of complexified zero-level sets of 
non-discriminant strictly Morse perturbations of a singularity, and the edges correspond to the flips of these data, which are defined by all possible surgeries of these perturbations. This graph can be created using a computer program. For the singularity $X_{10}^3$ (respectively, $X_{10}^1$), the graph has 3465452 (respectively, 1889028) vertices. 

We split this graph into subgraphs, whose vertices can be obtained from each other by the flips corresponding to the surgeries not related to the discriminant crossing.
These subgraphs are the ``virtual'' analogs of the connected components of the discriminant complement. In particular, if two real morsifications belong to the same component of the discriminant complement, then the corresponding vertices in the graph are in the same subgraph. In the present case of $X_{10}$ singularities, this part of the work was done partly by hand, but it can be fully automated.

One reason the the lists of subgraphs and of components of discriminant complements differ is that morsifications from different components can have the same set of topological data (``associated virtual morsifications''). This is usually associated with a symmetry group that acts on the argument space of our functions and preserves the original function like (\ref{meqx3}) or (\ref{meqx101}). The induced action of this group on the space of the versal deformation can permute the components of the discriminant complement while preserving the topological data. The resulting difference between the lists is expressed in the final statements of Theorems \ref{totalnum0} and \ref{mthm1010}: in our examples, there are 32 and 12 subgraphs, some of which correspond to pairs of symmetric components.

A more significant potential discrepancy between the lists stems from the fact that there is no guarantee that all standard flips of topological data can be always implemented by the real surgeries of functions. For this reasons, some vertices of the graph found by the algorithm may not correspond to any real morsification, and moreover there may be subgraphs all of whose vertices do not represent the real morsifications. However,
it turns out that  this inconsistency does not happen in the case of $X_{10}$ singularities. 
Indeed, we found a real morsification representing each subgraph obtained in this way, see \S \S \ref{real3} and \ref{real4}.

Another possible resulting problem is that the sets of topological data of morsifications which actually lie in different components of the discriminant complement can be connected by a chain of formal flips, and therefore belong to the same subgraph. Consequently, these real components will be counted as one. This is the most serious problem, it is the reason why we formulate Conjectures \ref{mcj1} and \ref{mcj2} and not theorems with the same assertions.

 These issues do not arise for simple singularities due to the {\em properness of the Lyashko--Looijenga map}, see \cite{Lo0}, nor for parabolic singularities, which is proven in \cite{vasX9}, \cite{vasP8}, \cite{vasJ10} based on P.~Jaworski's results (see \cite{Jaw} and \cite{Jaw2}). The $X_{10}$ singularities are the first example in  which analogous statements have  not been proven.

\subsection{Virtual morsifications and virtual surgeries}

Let $f:(\C^2, \R^2, 0) \to (\C, \R, 0)$ be a holomorphic function germ
having an isolated critical point at the origin and let $F:(\C^2 \times \C^k, \R^2 \times \R^k, 0) \to (\C, \R, 0)$ be its $k$-parametric versal deformation. The Milnor number $\mu(f)$ of this critical point is finite and is equal to the Milnor number of the function $f(x,y) + z^2: (\C^3,\R^3, 0) \to (\C, \R, 0)$.
Let $f_\lambda: (\C^2, \R^2) \to (\C, \R)$, $\lambda \in \R^k \setminus \Sigma(F),$ be a very small non-discriminant real perturbation of the function $f$. Denote by $V_\lambda$ the {\em Milnor fiber} of the perturbation $f_\lambda(x,y) + z^2$ of the function $f(x, y) + z^2$, i.e. the intersection of its zero-level set $(f_\lambda+ z^2)^{-1}(0) \subset \C^3$ with a small ball centered at the origin. It is known (see \cite{Milnor}, \cite{AVGZ}) that $H_2(V_\lambda) \simeq \Z^{\mu(f)}.$

\begin{figure}
\unitlength 1.1mm
\begin{center}
\begin{picture}(60,40)
\put(0,20){\line(1,0){60}}
\put(5,20){\circle*{1.2}}
\put(15,20){\circle*{1.2}}
\put(40,20){\circle*{1.2}}
\put(50,20){\circle*{1.2}}
\bezier{100}(30,20)(17,30)(5,20)
\bezier{50}(30,20)(22,24)(15,20)
\bezier{50}(30,20)(35,24)(40,20)
\bezier{100}(30,20)(40,30)(50,20)
\bezier{100}(30,20)(30,30)(35,35)
\bezier{100}(30,20)(30,10)(35,5)
\put(35,35){\circle*{1.2}}
\put(35,5){\circle*{1.2}}
\put(13,35){\circle*{1.2}}
\put(13,5){\circle*{1.2}}
\put(45,37){\circle*{1.2}}
\put(45,3){\circle*{1.2}}
\bezier{100}(30,20)(25,32)(13,35)
\bezier{100}(30,20)(25,8)(13,5)
\bezier{100}(30,20)(35,30)(45,37)
\bezier{100}(30,20)(35,10)(45,3)
\put(4,16){{\tiny 1}}
\put(14,16){{\tiny 2}}
\put(39,16){{\tiny 3}}
\put(49,16){{\tiny 4}}
\put(12,31){{\tiny 5}}
\put(32,35){{\tiny 6}}
\put(44,33){{\tiny 7}}
\put(12,7){{\tiny 8}}
\put(32,4){{\tiny 9}}
\put(44,5){{\tiny 10}}
\end{picture} \qquad \quad
\begin{picture}(60,40)
\put(0,20){\line(1,0){60}}
\put(5,20){\circle*{1.2}}
\put(15,20){\circle*{1.2}}
\put(40,20){\circle*{1.2}}
\put(50,20){\circle*{1.2}}
\bezier{100}(30,20)(17,30)(5,20)
\bezier{50}(30,20)(22,24)(15,20)
\bezier{50}(30,20)(35,24)(40,20)
\bezier{100}(30,20)(40,30)(50,20)
\bezier{100}(30,20)(30,30)(35,35)
\bezier{100}(30,20)(30,10)(35,5)
\put(35,35){\circle*{1.2}}
\put(35,5){\circle*{1.2}}
\put(13,35){\circle*{1.2}}
\put(13,5){\circle*{1.2}}
\put(45,37){\circle*{1.2}}
\put(45,3){\circle*{1.2}}
\bezier{100}(30,20)(25,32)(13,35)
\bezier{100}(30,20)(25,8)(13,5)
\bezier{100}(30,20)(25,30)(10,31)
\bezier{50}(10,31)(0,35)(10,37)
\bezier{100}(10,37)(35,37)(45,37)
\bezier{100}(30,20)(25,10)(10,9)
\bezier{50}(10,9)(0,5)(10,3)
\bezier{100}(10,3)(35,3)(45,3)
\put(4,16){{\tiny 1}}
\put(14,16){{\tiny 2}}
\put(39,16){{\tiny 3}}
\put(49,16){{\tiny 4}}
\put(12,31.3){{\tiny 6}}
\put(44,33){{\tiny 5}}
\put(35,32){{\tiny 7}}
\put(12,6.4){{\tiny 9}}
\put(44,5){{\tiny 8}}
\put(35,6){{\tiny 10}}
\end{picture}
\end{center}
\caption{Standard systems of paths}
\label{standd}
\end{figure}

Assume that the perturbation $f_\lambda$ is {\em strictly Morse} and non-discriminant, that is, it has exactly $\mu(f)$ critical points close to the origin in $\C^2$, all these points are Morse, and all corresponding critical values are distinct from each other and from 0. The set of these critical values is invariant under the complex conjugation in $\C^1$. Let us connect these values with the non-critical value 0 by a system of non-intersecting paths as shown in Fig.~\ref{standd}. The paths to the complex conjugate non-real critical values should be themselves complex conjugate, and the paths to the real values should go in the upper half-plane, staying sufficiently close to the real line. A basis of the homology group $H_2(V_\lambda)$ consists of {\em vanishing cycles} defined by these paths, see e.g. \cite{AVGZ}. If we fix an orientation of the space $\R^3$, then there is a canonical way to orient all of these basic cycles, see \S V.1.6 in \cite{ramif}. These cycles are also naturally ordered. In particular, the cycles that vanish at the real critical points come first in the ascending order by critical value.

\begin{definition} \rm
\label{virt}
The {\it virtual morsification associated} with the strictly Morse non-discriminant perturbation $f_\lambda$ is the collection of its discrete data consisting of:
\begin{enumerate}
\item[(a)] the $\mu(f) \times \mu(f)$ matrix of intersection indices of basic vanishing cycles in the complex manifold $V_\lambda$,
\item[(b)] the string of intersection indices of these vanishing cycles with the naturally oriented {\em real cycle} $V_\lambda \cap \R^3$,
\item[(c)] the string of {\em positive} Morse indices (i.e., the numbers of positive squares in the Morse normal form)
of the real critical points of $f_\lambda$ at which these cycles vanish, and
\item[(d)] the numbers of negative and positive critical values of $f_\lambda$.
\end{enumerate}
\end{definition}

For examples of virtual morsifications written in the canonical form, see (\ref{7200}) 
\begin{equation} 
\begin{array}{|ccccc|ccccc||}
\hline
 $-2$ & 0 & 0 & 0 & $-2$ & 1 & 1 & 0 & 1 & 0 \\
 0 & $-2$ & 0 & 0 & 0 & 0 & 0 & 1 & 0 & 0 \\
 0 & 0 & $-2$ & 0 & 0 & 1 & 0 & 0 & 0 & 1 \\ 
 0 & 0 & 0 & $-2$ & 0 & 0 & 1 & 1 & 0 & 0 \\
 $-2$ & 0 & 0 & 0 & $-2$ & 1 & 1 & 0 & 1 & 0 \\
 1 & 0 & 1 & 0 & 1 & $-2$ & 0 & 0 & 0 & 0 \\
 1 & 0 & 0 & 1 & 1 & 0 & $-2$ & 0 & 0 & 0 \\
 0 & 1 & 0 & 1 & 0 & 0 & 0 & $-2$ & 0 & 0 \\
 1 & 0 & 0 & 0 & 1 & 0 & 0 & 0 & $ -2$ & 0 \\
 0 & 0 & 1 & 0 & 0 & 0 & 0 & 0 & 0 & $-2$ \\
\hline
\hline
2 & 2 & 2 & 2 & 2 & $-2$ & $-2$ & $-2$ & $-2$ & $-2$ \\
\hline
\hline
1 & 2 & 2 & 2 & 2 & 1 & 1 & 1 & 1 & 1 \\
\hline
\end{array}
\label{7200}
\end{equation}
(which is associated with the real morsification shown in the picture $(7200, \ 3)$)
and also formulas (\ref{m4800}) and (\ref{m4320}) on pages \pageref{m4800} and \pageref{m4320} respectively.
The vertical lines in these tables separate the data for cycles vanishing in the points with negative, positive and non-real critical values. For example, the element (d) of the virtual morsification (\ref{7200}) (i.e., the numbers of negative and positive critical values) is equal to $(5, 5).$ In all three of our examples, all critical points are real. This fact is indicated by a vertical line in the right margin of the table: it separates  the empty set of imaginary critical points. 
If not all critical points of $f_\lambda$ are real, then the bottom row of the table (expressing the element (c) of the virtual morsification) will have several blank places on the right. The numbers 1 and 2 on the bottom line indicate, respectively, the saddlepoints and the local minima.

\begin{remark} \rm
If the number of pairs of imaginary critical points is greater than one, then there can be more than one virtual morsification associated with $f_\lambda$, because in this case there is no canonical way to draw the paths to the imaginary critical values: see two parts of Fig.~\ref{standd}. 

If all $\mu(f)$ critical points of $f_\lambda$ are real, then element (b) of the associated virtual morsification is determined by its other elements, see formulas (V.6), (V.8) of \cite{ramif}.
\end{remark}

\begin{definition} \rm
\label{defsur}
{\em Elementary surgeries} of virtual morsifications include six flips of these data reflecting the topological surgeries of the corresponding real morsifications, namely
\begin{itemize}
\item[(s1, s2)] collision of two neighboring real critical values at a value different from 0, after which the corresponding two critical points either (s1) meet at a degenerate critical point and go into the complex domain or (s2) change the order  of their critical values in $\R^1$, 
\item[(s3, s4)] collision of two complex conjugate critical values at a point on the line $\R^1 \setminus \{0\}$, after which the corresponding critical points of $f_\lambda$ in $\C^2$ either (s3) meet at a real point and enter $\R^2$ or (s4) miss one another while the imaginary parts of their critical values change their signs, 
\item[(s5, s6)] jump of a real critical value up (s5) or down (s6) through 0; \\ \hspace*{4cm} and additionally
\item[(s7)] specifically virtual surgeries within the classes of virtual morsifications associated with the same real morsification, caused by changes of systems of paths going from 0 to imaginary critical values, see two parts of Fig.~\ref{standd}.
\end{itemize}
\end{definition}

The information contained in the virtual morsification associated with any real morsification $f_\lambda$ is sufficient to predict the virtual morsification associated with the real morsification obtained from $f_\lambda$ by arbitrary surgeries of types (s1)--(s6) that are not prohibited by this information. See \S V.9 in \cite{ramif} for explicit formulas relating the original virtual morsification to the new one. For example, the numbers contained in the (d) elements of the virtual morsifications are unchanged by surgeries (s2), (s4), and (s7). They change by $\pm 1$ after the surgeries (s5/s6), by 0 and $- 2$ after the surgery (s1) and by 0 and +2 after the surgery (s3). However, this information generally does not allow one to determine if a particular type of real surgery can be performed: it only provides the necessary conditions for it. For simple and parabolic singularities, these conditions are also sufficient, but for more complicated singularities this is generally not the case.

\begin{definition} \rm
A {\em virtual morsification} of the singularity $f$ is any collection of data as in Definition \ref{virt}, that can be obtained from a virtual morsification associated with an arbitrary non-discriminant strict morsification of $f$ by a chain of formal flips modeling the elementary surgeries (s1)--(s7).

These flips are called the {\em virtual surgeries} of virtual morsifications.

The {\it formal graph} of a singularity $f$ is a graph whose vertices are all the virtual morsifications of $f$, and two vertices are connected by an edge if and only if the corresponding virtual morsifications can be obtained from each other by a single virtual surgery.

A {\em virtual component} of the formal graph is an arbitrary its subgraph obtained from it by removing all edges of type (s5/s6), i.e., the virtual surgeries that model the surgeries of $f_\lambda$ related to crossing the discriminant variety.
\end{definition}

\begin{proposition}
\label{taut}
The formal graph of a singularity is independent of the choice of the initial strict morsification. 
All virtual morsifications associated with real morsifications of a function singularity appear as vertices in its formal graph. 
If two parameter values $\lambda$ and $\tilde \lambda$ of a versal deformation
 belong to the same connected component of the complement of the real discriminant then all the virtual morsifications associated with the perturbations $f_\lambda$ and $f_{\tilde \lambda}$ belong to the same virtual component.
\end{proposition}

\noindent
{\it Proof.} The parameter space of any versal deformation is a connected piece of a vector space. The set of discriminant or non-Morse functions is a semialgebraic subvariety of codimension one in the parameter space. Therefore, any two generic morsifications of the same singularity can be connected by a path in this space that crosses these bifurcation sets only finitely many times. The corresponding vertices of the formal graph are then connected by a sequence of edges corresponding to these bifurcations.

If the morsifications belong to the same component of the discriminant complement, they can be connected by a path with no edges of types s5 or s6. \hfill $\Box$

\subsection{Basic invariants of real and virtual morsifications}
\label{basinv} 

\begin{definition} \rm
The invariant \ $\mbox{Ind}$ of a virtual morsification is equal to the number of even indices in the left part of the bottom row of its canonical table like (\ref{7200}) minus the number of odd indices in the same part of the row.
\end{definition}

For example, this part of the virtual morsification (\ref{7200}) consists of one odd number and four even numbers; therefore, the value of $\mbox{Ind}$ is equal to 3.

\begin{proposition}
The invariant $\mbox{Ind}$ of a real morsification $($see Definition \ref{defind}$)$ is equal to the  invariant $\mbox{Ind}$ of any associated virtual morsification.
\end{proposition}

\noindent
{\it Proof} follows immediately from the definitions. \hfill $\Box$
 
\begin{definition}
\rm The {\em set-valued invariant} 
of a component of the complement of the real discriminant is the set of vertices of the virtual component containing the virtual morsifications associated with the real morsifications in this component. The invariant \ $\mbox{Card}$ \ of such a component is the number of elements in this set.
\end{definition}

\subsection{The algorithm}

Our investigation of the discriminant complement of a singularity consists of three steps:

1. Compute the formal graph of the singularity and enumerate all its virtual components.

2. Choose one virtual morsification from each virtual component and try to find a real morsification with these topological data. According to Proposition \ref{taut} all obtained real morsifications will represent different components or $\R^k \setminus \Sigma(F)$;

3. Think whether the obtained list of real morsifications represents all components of the real discriminant complement.
\medskip

The first step is purely combinatorial and can, in principle, be solved by a computer. 

The second step of the algorithm is much less formal. However, in all the examples considered thus fr, it has been possible to accomplish it, see \cite{ramif}, \cite{VLoo}, \cite{parab}, and \S\S \ref{real3}, \ref{real4} of the present paper. 

The question contained in the third step has a positive answer for simple and parabolic singularities; however,   the answer is unknown  for $X_{10}$.

\begin{remark} \rm In the future, more complicated examples may present a new difficulty can appear: the formal graph of the singularity may be infinite.
In this case, we can define an increasing sequence of finite subgraphs such that each vertex of the whole formal graph appears in a subgraph of the sequence. For example, for each natural number $N$, we can define the set of vertices of the $N$th subgraph as the set of virtual morsifications, all of whose numerical elements are in the interval $[-N, N]$. Then solve the problem for several initial subgraphs of this sequence and analyze the stabilization of results.
\end{remark}

Information on formal graphs and their subgraphs for $X_{10}$ singularities is summarized in the following two theorems.

\begin{table}
\caption{Numbers of virtual morsifications for the $X_{10}^3$ singularity}
\label{taX103}
{\footnotesize
\begin{tabular}{|c|c|c|c|c|c|c|c|c|c||c|}
\hline
Ind & $-5$ & $-4$ & $-3$ & $-2$ & $-1$ & 0 & 1 & 2 & 3 & Total \\
\hline
$\#$ & 18192 & 158792 & 457292 & 779716 & 529728 & 712396 & 493384 & 255520 & 60432 & 3465452 \\
\hline
\end{tabular}
}
\end{table}

\begin{theorem}
\label{totalnum}
1. The singularity $X_{10}^3$ has exactly 3465452 distinct virtual morsifications. The invariant \ $\mbox{Ind}$ \ takes on them values from \ $-5$ \ to \ $3$. The numbers of virtual morsifications with different values of this invariant are shown in Table \ref{taX103}. 

2. These virtual morsifications are divided into exactly 32 virtual components.
The values of invariants \ $\mbox{Card}$ \ and \ $\mbox{Ind}$ \ for these components are indicated in the subscripts of the 32 pictures in Figs.~\ref{x3beg}--\ref{x3end}. 

3. For each of these 32 pairs of basic invariants, the corresponding virtual component contains a virtual morsification associated with the real morsification shown in the picture above this pair of indices. 
\end{theorem}


{\rm \begin{table}
\caption{Numbers of virtual morsifications for the $X_{10}^1$ singularity}
\label{taX101}
{\footnotesize
\begin{tabular}{|c||c|c|c|c|c|c|c|c|c||c|}
\hline
$|\mbox{Ind}|$ & $-4$ & $-3$ & $-2$ & $-1$ & 0 & 1 & 2 & 3 & 4 & Total \\
\hline
$\#$ & 40032 & 157004 & 265228 & 306917 & 541540 & 277887 & 202812 & 84936 & 12672 & 1889028 \\
\hline
\end{tabular}
}
\end{table}
}

\begin{theorem}
\label{mthmx101}
1. The formal graph of the singularity $X_{10}^1$ contains 1889028 distinct virtual morsifications. The invariant \ {\rm $\mbox{Ind}$} \ takes values from $-4$ to $4$ on its elements. The numbers of virtual morsifications with these values of \ {\rm $\mbox{Ind}$ } \ are shown in Table~\ref{taX101}.

2. There are exactly 12 distinct virtual components of the formal graph of the singularity $X_{10}^1$. The basic invariants $\mbox{Card}$ and $\mbox{Ind}$ of these virtual components are indicated under the pictures of Fig. \ref{tta}.

3. The 12 virtual components are represented by 12 virtual morsifications, which are associated with real morsifications whose zero-level sets are shown in the pictures of Fig. \ref{tta} above the corresponding values of the basic invariants.
\end{theorem}

\section{Proof of Theorems \ref{totalnum} and \ref{mthmx101}}

\subsection{Scheme of the proof}

Statements 1 and 2 of these theorems are derived from computer calculations. To begin, we take an arbitrary sabirization $f$ of the $X_{10}^3$ or $X_{10}^1$ singularity, see Definition \ref{sabdef}. Using the method of \cite{AC}, \cite{GZ}, we calculate the intersection matrix of the vanishing cycles of a non-discriminant strictly Morse perturbation $\tilde f(x, y) + z^2$ of the function $f(x, y) + z^2$. Using formulas (V.6) and (V.8) from \cite{ramif}, we deduce the intersection indices of these vanishing cycles with the real cycle $V_\lambda \cap \R^3$, i.e., the element (b) of the virtual morsification associated with $\tilde f$. We then substitute this virtual morsification as the initial data for a combinatorial program that counts all other virtual morsifications by performing arbitrary sequences of virtual surgeries. This program, with different initial data, is available in \cite{pro2}, and its detailed description an be found in \S \S V.8--9 of \cite{ramif} or in \S 5.3 of \cite{AGLV}. In particular, it counts the number of virtual morsifications with any value of the invariant $\mbox{Ind}$. The results of the computations are presented in Tables \ref{taX103} and \ref{taX101} and prove the first statements of our theorems. Next, we use a {\em reduced version} of the program that does not allow surgeries of types s5 and s6. Starting from an arbitrary virtual morsification, the reduced program counts all virtual morsifications in its virtual component and gives the value of the invariant $\mbox{Card}$ for the elements of this component. Next, we find a virtual morsification that is not in this virtual component and repeat the process until the sum of the values of $\mbox{Card}$ over all found virtual components equals the total number of virtual morsifications found in statement 1. The results of this procedure prove the second statements of the theorems.

The proof of the third statements consists of two steps:

a) realizing all 32 pictures of Figs.~\ref{x3beg}--\ref{x3end} and 12 pictures of Fig.~\ref{tta} by real morsifications and 

b) proving that the subscript of each of these pictures contains the values of the basic invariants of the morsifications shown in this picture. 

Step a) is performed in \S \ref{real3} for the case $X_{10}^3$ and in \S \ref{real4} for the case $X_{10}^1$. The part of step b) concerning the $\mbox{Ind}$ invariant follows immediately from these realizations, the part concerning the $\mbox{Card}$ invariant occupies the rest of the paper.

\subsection{Virtual morsifications and $\mbox{Card}$ invariants associated with real morsifications of the $X_{10}^3$ and $X_{10}^1$ singularities}

\label{realform}

\subsubsection{} \label{sabre}
Gusein-Zade--A'Campo method \cite{AC}, \cite{GZ} allows us to calculate the virtual morsifications associated with all 25 perturbations of the $X_{10}^3$ singularity constructed in \S \ref{realsab} and with all 12 perturbations of the $X_{10}^1$ singularity constructed in \S \ref{real4}. The reduced program applied to all these virtual morsifications proves that
the first subscripts under the corresponding pictures in
Figs.~\ref{x3beg}--\ref{x3end} and \ref{tta} are indeed equal to the values of 
 the invariant \ $\mbox{Card}$ of these perturbations. In the case of the $X_{10}^1$ singularity, these 12 perturbations represent all of the found components; therefore, Theorem \ref{mthmx101} is proved.
It remains to prove the same for seven perturbations of $X_{10}^3$ singularity that are only realized in \S \ref{realadj}.

\subsubsection{Virtualization of the morsification shown in picture $(7200, \ 3)$}
\label{ssreal7200}

According to statement 2 of Theorem \ref{totalnum}, there are only two virtual components with $\mbox{Ind}=3$, namely $(7200, \ 3)$ and $(53232, \ 3)$. The representatives of all virtual components are obtained when proving this statement 2. 

 Upon a special request,  the reduced program checks that the virtual component (53232, \ 3) does not contain virtual morsifications, whose critical point with the lowest critical value is a saddlepoint. The real morsification constructed in \S \ref{adje8} and shown in the picture (7200, \ 3) has such a critical point. Therefore, the virtual morsification associated with this real morsification belongs to the only other virtual component $(7200, \ 3)$.

\subsubsection{Virtualization of picture $(4800, \ -5)$}

\label{ss4800}

According to statement 2 of Theorem \ref{totalnum}, there are only two virtual components with $\mbox{Ind} = -5$. Their Card invariants are  4800 and 13392.

The morsification constructed in \S \ref{adje8} and realizing the picture $(4800, \ -5)$ has $\mbox{Ind}=-5$ and satisfies the following conditions:

\begin{itemize}
\item
It has ten real critical points, none of which are local minima.
\item
It has exactly five critical points with negative critical values.
\item
Its critical points with the greatest and the smallest critical values are saddlepoints.
\end{itemize}

Our reduced program verifies that the virtual component with $\mbox{Card} = 13392$ does not contain virtual morsifications with these properties.
Therefore, this real morsification has $\mbox{Card}=4800$.
\medskip

The virtual component $(4800, \ -5)$ of singularity $X_{10}^3$ contains the following virtual morsification: 
\begin{equation} 
\begin{array}{|ccccc|ccccc||}
\hline
 $-2$ & 0 & 0 & 0 & 0 & 0 & 1 & 0 & 0 & 1 \\
 0 & $-2$ & 0 & 0 & 0 & 0 & 0 & 1 & 1 & 1 \\
 0 & 0 & $-2$ & 0 & 0 & 1 & 1 & 0 & 0 & 1 \\
 0 & 0 & 0 & $-2$ & 0 & 0 & 1 & 1 & 0 & 2 \\
 0 & 0 & 0 & 0 & $-2$ & 0 & 0 & 0 & 1 & 0 \\
 0 & 0 & 1 & 0 & 0 & $-2$ & 0 & 0 & 0 & 0 \\
 1 & 0 & 1 & 1 & 0 & 0 & $-2$ & 0 & 0 & $-2$ \\
 0 & 1 & 0 & 1 & 0 & 0 & 0 & $-2$ & 0 & $-2$ \\
 0 & 1 & 0 & 0 & 1 & 0 & 0 & 0 & $-2$ & 0 \\
 1 & 1 & 1 & 2 & 0 & 0 & $-2$ & $-2$ & 0 & $-2$ \\
\hline
\hline
0 & 0 & 0 & 0 & 0 & 0 & 0 & 0 & 0 & -6 \\
\hline
\hline
1 & 1 & 1 & 1 & 1 & 0 & 0 & 0 &0 & 1 \\
\hline
\end{array}
 \ . \label{m4800}
\end{equation}

\begin{lemma}
\label{lell}
All virtual morsifications in the virtual component $(4800, \ -5)$ can be obtained from $($\ref{m4800}$)$ by  surgeries of only type s2. 
\end{lemma}

\noindent
{\it Proof.} Real surgery s2 permutes two critical points with neighboring critical values of the same sign if the intersection index of the corresponding vanishing cycles is equal to 0.
Corresponding virtual surgery permutes the two corresponding lines and rows of the intersection matrix. Applying these operations to matrix (\ref{m4800}) yields 4800 different virtual morsifications. Indeed, the critical points with negative critical values can be permuted  arbitrarily, this gives the factor $5!=120$ of the number of all admissible permutations. The first four critical points with positive critical values can be permuted in any way; this yields the summand $4!=24$ of the second factor. In addition, the critical value of any of the critical points No. 6, No. 9, or both, can be moved above the critical value of critical point 10. In these three cases, there are 6, 6, and 4 choices, respectively, of ordering the maximum points. In total, this yields $120 \times (24 + 6 + 6 + 4) = 4800$ admissible permutations.

The Dynkin graph of the quadratic form defined by matrix (\ref{m4800}) has no symmetries; therefore, all the described permutations define different virtual morsifications and hence exhaust the entire virtual component. \hfill $\Box$ 

\begin{corollary}
\label{cor87}
1. For each strict morsification in the component $(4800, \ -5)$, the two smallest positive critical values are attained at the points of local maxima.

2. All real morsifications in this component have real critical points only. \hfill $\Box$
\end{corollary}

\subsubsection{Virtualization of picture $(4320, \ -4)$}
\label{vrt4320}

According to statement 2 of Theorem \ref{totalnum}, there are only six virtual components with $\mbox{Ind} = -4$. Their Card invariants are 8496, 9984, 25104, 81848, 29040, and 4320.

The morsification constructed in \S \ref{tox92} and realizing the picture $(4320, \ -4)$
satisfies the following conditions:

\begin{itemize}
\item
It has ten real critical points, one of which is a local minimum and three are local maxima.
\item
It takes its smallest critical value at the minimum point, and takes the seventh, eighth, and ninth critical values in ascending order at the maxima.
\item
It has exactly six critical points with negative critical values.
\item
The intersection indices of nine vanishing cycles corresponding to all critical points except for the minimum point and ordered by the ascending of critical values form the matrix (\ref{4320a}).
\end{itemize}

Our reduced program verifies that the virtual components with $\mbox{Card}$ equal to 81848, 25104, 29040, 9984, or 8496 do not contain virtual morsifications with these properties. Therefore, the $\mbox{Card}$ invariant of the morsification shown in picture $(4320, \ -4)$ is indeed 4320. 

The virtual component $(4320, \ -4)$ contains the following virtual morsification: 
\begin{equation} 
\begin{array}{|cccccc|cccc||}
\hline
 $-2$ & 0 & 0 & 0 & 0 & 0 & 0 & 0 & 0 & 1 \\
 0 & $-2$ & 0 & 0 & 0 & 0 & 0 & 0 & 1 & 0 \\
 0 & 0 & $-2$ & 0 & 0 & 0 & 0 & 1 & 1 & 1 \\
 0 & 0 & 0 & $-2$ & 0 & 0 & 0 & 1 & 0 & 1 \\
 0 & 0 & 0 & 0 & $-2$ & 0 & 1 & 1 & 0 & 1 \\
 0 & 0 & 0 & 0 & 0 & $-2$ & 1 & 0 & 0 & 0 \\
 0 & 0 & 0 & 0 & 1 & 1 & $-2$ & 0 & 0 & 0 \\
 0 & 0 & 1 & 1 & 1 & 0 & 0 & $-2$ & 0 & $-2$ \\
 0 & 1 & 1 & 0 & 0 & 0 & 0 & 0 & $-2$ & 0 \\
 1 & 0 & 1 & 1 & 1 & 0 & 0 & $-2$ & 0 & $-2$ \\
\hline
\hline
2 & 0 & 0 & 0 & 0 & 0 & 0 & 0 & 0 & $-4$ \\
\hline
\hline
2 & 1 & 1 & 1 & 1 & 1 & 0 & 0 & 0 & 1 \\
\hline
\end{array} \ . \label{m4320}
\end{equation}

\begin{lemma}
\label{lellu}
All virtual morsifications in the virtual component $(4320, \ -4)$ can be obtained from $($\ref{m4320}$)$ by surgeries of only type s2. 
\end{lemma}

\noindent
{\it Proof.} By applying these surgeries to the intersection matrix (\ref{m4320}), we can obtain 4320 different virtual morsifications. Indeed, the compositions of real s2 surgeries can arbitrarily permute the critical points with negative critical values. This gives the factor of $6!=720$ of the number of admissible permutations. The first three critical points with positive critical values can be permuted in any way, this yields a summand of $3!=6$ of the second factor. Additionally, the critical value of  critical points No. 6, or No. 9, or both of them can be moved above the value of the critical point 10; in each of these three cases, there are two admissible choices for ordering the maximum points. The resulting number of permutations, $720 \times (6 + 2 + 2 + 2)= 8640$,  should be divided by two, since the Dynkin graph of the quadratic form (\ref{m4320}) has a symmetry of order two; thus, each permutation thus obtained has exactly one other permutation that defines the same virtual morsification. \hfill $\Box$ 

\begin{corollary}
\label{cor87z}
1. For each strict morsification in the component $(4320, \ -4)$, the smallest positive critical value is attained at a point of local maximum.

2. All real morsifications in the component with $\mbox{Card}=4320$ have only real critical points. 
\hfill $\Box$
\end{corollary}
 
\subsubsection{The virtualization of pictures $(25104, \ -4)$ and $(63700, -3)$}

The morsification shown in picture $(25104, \ -4)$ can be obtained from a morsification $\tilde f$ shown in the picture $(4800, \ -5)$ by a single surgery of type s6, which consists in adding a constant function. In this surgery, the smallest positive critical value passes the zero level and becomes the greatest negative critical value of the resulting function. Virtual morsifications obtained from each other by the s5/s6 surgery can easily be computed one through the other: the new matrix of intersection indices is obtained from the old one by  transposing the neighboring rows and columns. In \S~\ref{ss4800}, we proved that the virtual morsification associated with the real morsification $\tilde f$ belongs to the virtual component $(4800, \ -5)$. If this virtual morsification is equal to (\ref{m4800}), then the reduced program applied to the resulting virtual morsification shows that its Card invariant is indeed equal to 25104. 

The same result occurs if we 
modify the original virtual morsification (\ref{m4800}) inside the virtual component $(4800, \ -5)$ by arbitrary s2 surgeries that permute critical points with negative critical values, and also by the allowed s2 surgeries that permute critical points with the four highest critical values. Indeed, all these modifications commute with the s6 surgery. The only modifications to the original virtual morsification that could change the resulting virtual component are those that change the critical point with the smallest positive critical value (which is intended to pass the zero level). According to Corollary \ref{cor87}, this smallest positive critical value is attained at a maximum point. 
The reduced program examined the resulting virtual morsifications for all three possible choices of the maximum point made the lowest by these modifications 
and verified that all of them belong to the virtual component with Card = 25104. Thus, the real morsification shown in picture $(25104, \ -4)$ has this value of the Card invariant. 

Analogously, the morsification shown in picture $(63700, \ -3)$ can be obtained from a morsification from the real component $(4800, \ -5)$ through two surgeries s6, which change the signs of the two lowest positive critical values (which, according to Corollary \ref{cor87}, are attained at local maxima). The virtual component of the virtual morsification associated with the obtained real morsification depends only on the choice of the two corresponding  critical points. The reduced program considered all six choices of these points and verified that, regardless of the choice, the resulting virtual morsification belongs to the virtual component with Card = 65700. Thus, the real morsification shown in picture $(63700, \ -3)$ has the Card invariant 65700. 

\subsubsection{Virtualization of the picture $(22800, \ -3)$}

Each morsification shown in picture $(22800, \ -3)$ can be obtained from the morsification $\hat f$ shown in the picture $(4320, \ -4)$ through a single surgery of type s6. In this surgery, the smallest positive critical value passes the zero value and becomes the greatest negative critical value of the resulting function. According to \S~\ref{vrt4320}, the virtual morsification associated with $\hat f$ belongs to the virtual component $(4320, \ -4)$. The reduced program shows that if this virtual morsification is equal to (\ref{m4320}), then the virtual morsification obtained from it by the virtual s6 surgery has $\mbox{Card}=22800$. Hence, the the same is true for the associated real morsification shown in picture $(22800, \ -3)$. 

This conclusion holds if we modify the original virtual morsification inside the virtual component $(4320, \ -4)$ by arbitrary s2 surgeries that permute the negative critical values, and also by the allowed s2 surgeries that permute the upper three critical values. Indeed, each such modification commutes with the virtual s6 surgery. The only modifications to the original virtual morsifications that could alter the resulting virtual component are those that alter the critical point with the smallest positive critical value (which passes through zero). 
The reduced program examined the obtained virtual morsifications for all two possible choices of this new critical point  and verified that they still belong to the virtual component with Card = 22800 in all cases. Thus, the real morsification shown in picture $(22800, \ -3)$ has this value of Card invariant. 

\subsubsection{Virtualization of the picture $(8496, \ -4)$ }
\label{vrt8496}

 The real morsification shown in the picture $(8496, \ -4)$ has $\mbox{Ind} = -4$. Therefore, by statement 2 of Theorem \ref{totalnum}, its $\mbox{Card}$ invariant can be equal only to one of six values 8496, 9984, 25104, 81848, 29040, or 4320. 
 
\begin{figure}
\unitlength 1mm
\begin{picture}(120,120)
\put(0,0){\line(1,0){120}}
\put(0,0){\line(0,1){120}}
\put(120,120){\line(-1,0){120}}
\put(120,120){\line(0,-1){120}}
\thicklines
\bezier{2000}(20,120)(70,65)(30,80)
\bezier{1500}(30,80)(5,90)(5,60)
\bezier{1500}(5,60)(5,30)(30,40)
\bezier{2000}(30,40)(70,55)(20,0)
\bezier{3000}(120,10)(100,60)(120,110)
\bezier{4000}(110,0)(70,60)(110,120)
\put(73,62){\makebox(0,0)[cc]{\tiny $+$}}
\put(63,42){\makebox(0,0)[cc]{\tiny $+$}}
\put(63,78){\makebox(0,0)[cc]{\tiny $+$}}
\put(71.7,60){\circle{8.32}}
\put(65,40){\circle{8.32}}
\put(65,80){\circle{8.32}}
\put(71.7,60){\circle{8.35}}
\put(65,40){\circle{8.35}}
\put(65,80){\circle{8.35}}
\put(71.7,60){\circle{8.33}}
\put(65,40){\circle{8.33}}
\put(65,80){\circle{8.33}}
\thinlines
\put(56.5,80){\vector(-1,0){8.5}}
\put(57,80){\vector(1,0){8.5}}
\put(56.5,40){\vector(-1,0){8.5}}
\put(57,40){\vector(1,0){8.5}}
\put(56.5, 80){\circle*{1}}
\put(56.5, 40){\circle*{1}}
\put(68.5,50){\vector(1,3){3.3}}
\put(68.5,50){\vector(-1,-3){3.3}}
\put(68.5,70){\vector(-1,3){3.3}}
\put(68.5,70){\vector(1,-3){3.3}}
\put(68.5,70.5){\circle*{1}}
\put(68.5,49.5){\circle*{1}}
\put(25,60){\circle*{1}}
\bezier{150}(100,60)(100,90)(115,115)
\put(115,115){\vector(2,3){3}}
\bezier{150}(100,60)(100,30)(115,5)
\put(115,5){\vector(2,-3){3}}
\put(100,60){\circle*{1}}
\put(120,60){\vector(-1,0){20}}
\bezier{100}(25,60)(55,67)(56.5,80)
\put(55.5,76.5){\vector(1,3){1}}
\bezier{100}(25,60)(55,53)(56.5,40)
\put(55.5,43.5){\vector(1,-3){1}}
\bezier{150}(25,60)(50,57)(67,50)
\put(66,50.5){\vector(3,-1){2}}
\bezier{150}(25,60)(50,63)(67,70)
\put(66,69.5){\vector(3,1){2}}
\put(80,60){\vector(-1,0){8}}
\put(80,60){\circle*{1}}
\bezier{200}(40,120)(55,100)(56.5,81)
\put(56.5, 82.5){\vector(0,-1){2}}
\bezier{200}(40,0)(55,20)(56.5, 39)
\put(56.5,37.5){\vector(0,1){2}}
\bezier{250}(70,120)(80,75)(68.5,70.3)
\put(70,71){\vector(-2,-1){1}}
\bezier{250}(70,0)(80,45)(68.5,49.7)
\put(70,49){\vector(-2,1){1}}
\bezier{250}(90,120)(81,90)(80,61)
\put(80,62){\vector(0,-1){2}}
\bezier{250}(90,0)(81,30)(80,59)
\put(80,58){\vector(0,1){2}}
\bezier{150}(80,60)(97,60)(99,70)
\bezier{250}(99,70)(100,94)(117,120)
\put(116.5,118.5){\vector(1,3){0.4}}
\put(99,70){\vector(0,1){2}}
\bezier{250}(91,0)(82,30)(82,50)
\bezier{150}(100,60)(81,60)(82,50)
\put(98,60){\vector(1,0){1}}
\put(25.5,57.5){\makebox(0,0)[cc]{\tiny 1}}
\put(101.5,58){\makebox(0,0)[cc]{\tiny 10}}
\put(82,61.3){\makebox(0,0)[cc]{\tiny 2}}
\put(71,69){\makebox(0,0)[cc]{\tiny 3}}
\put(54,41.7){\makebox(0,0)[cc]{\tiny 4}}
\put(55,81.7){\makebox(0,0)[cc]{\tiny 5}}
\put(71,50.3){\makebox(0,0)[cc]{\tiny 6}}
\put(65,38){\makebox(0,0)[cc]{\tiny 7}}
\put(65,82){\makebox(0,0)[cc]{\tiny 8}}
\put(69,60){\makebox(0,0)[cc]{\tiny 9}}
\put(115,70){\small $-$}
\put(105,80){\small $+$}
\put(57,115){\small $-$}
\put(12,59){\small $-$}
\put(3,115){\small $+$}
\end{picture}
\caption{Separatrice diagram for perturbation (8496, \ $-4$)}
\label{separ}
\end{figure}

By the construction of this morsification given in \S \ref{tox91} (see especially Fig.~\ref{x91pert}), the diagram of separatrices of its gradient vector field in a generic metric has the structure shown in Fig.~\ref{separ}. In this picture, thick curves without arrows represent the zero-level set of the function. The thin curves with arrows are separatrices. The numbers at the critical points indicate the ascending order of the corresponding critical values. The way in which the separatrice going from the point 2 to the right and the separatrice coming to the point 10 from the left miss one another can be different depending on the choice of the generic metric in $\R^2$.

The associated virtual morsification satisfies the following properties: 

\begin{itemize}
\item
It has six saddlepoints, one local minimum, and three maxima.
\item
The critical point with the lowest critical value is a minimum. The maximum points have the critical values No. 7, 8, and 9.
\item 
It has exactly four critical points with positive critical values.
\item
The intersection indices of vanishing cycles corresponding to the lowest nine critical values are given by the matrix
{\rm 
\begin{equation} 
\begin{array}{ccccccccc}
 $-2$ & 0 & 1 & 1 & 1 & 1 & $-1$ & $-1$ & $-1$ \\
 0 & $-2$ & 0 & 0 & 0 & 0 & 0 & 0 & 1 \\ 
 1 & 0 & $-2$ & 0 & 0 & 0 & 0 & 1 & 1 \\
 1 & 0 & 0 & $-2$ & 0 & 0 & 1 & 0 & 0 \\
 1 & 0 & 0 & 0 & $-2$ & 0 & 0 & 1 & 0 \\
 1 & 0 & 0 & 0 & 0 & $-2$ & 1 & 0 & 1 \\
 $-1$ & 0 & 0 & 1 & 0 & 1 & $ -2$ & 0 & 0 \\
 $-1$ & 0 & 1 & 0 & 1 & 0 & 0 & $ -2$ & 0 \\
 $-1$ & 1 & 1 & 0 & 0 & 1 & 0 & 0 & $-2$ \\
\end{array}
 \ . \label{8496}
\end{equation} }
\end{itemize}

Indeed, this intersection matrix was computed by the Gusein-Zade--A'Campo method applied to the perturbation of $X_{9}^1$ singularity shown in Fig.~\ref{x91pert} (bottom).

The reduced program verifies that none of the virtual morsifications within the virtual components with $\mbox{Ind}=-4$, except for $(8496, \ -4),$ satisfy these conditions.
Thus, Theorem \ref{totalnum} is completely proved.

\medskip

\centerline{\bf Acknowledgment}
\smallskip

I am very grateful to the Weizmann Institute, where this work was carried out, and especially to Sergey Yakovenko for the extremely timely help and hospitality.

\end{document}